\documentclass[11pt]{article} 
\usepackage{amsmath}
\usepackage{amsfonts}
\usepackage{ifthen}
\usepackage{graphicx}
\usepackage{color}
\usepackage[latin1]{inputenc}
\usepackage[T1]{fontenc}
\usepackage[french,english]{babel}
\makeatletter
\def\@sect#1#2#3#4#5#6[#7]#8{%
  \ifnum #2>\c@secnumdepth
    \let\@svsec\@empty
  \else
    \refstepcounter{#1}%
    \protected@edef\@svsec{\@seccntformat{#1}\relax}%
  \fi
  \@tempskipa #5\relax
  \ifdim \@tempskipa>\z@
    \begingroup
      #6{%
        \@hangfrom{\hskip #3\relax\@svsec}%
          \interlinepenalty \@M #8\@@par}%
    \endgroup
    \csname #1mark\endcsname{#7}%
    \addcontentsline{toc}{#1}{%
      \ifnum #2>\c@secnumdepth \else
        \protect\numberline{\csname the#1\endcsname.}%
      \fi
      #7}%
  \else
    \def\@svsechd{%
      #6{\hskip #3\relax
      \@svsec #8}%
      \csname #1mark\endcsname{#7}%
      \addcontentsline{toc}{#1}{%
        \ifnum #2>\c@secnumdepth \else
          \protect\numberline{\csname the#1\endcsname.}%
        \fi
        #7}}%
  \fi
  \@xsect{#5}}
\def\@seccntformat#1{\csname the#1\endcsname.\quad}
\makeatother
\newtheorem{theo}[equation]{Th\'eor\`eme}
\newtheorem{etheo}[equation]{Theorem}
\newtheorem{proposition}[equation]{Proposition}
\newtheorem{fait}[equation]{Fait}
\newtheorem{lem}[equation]{Lemme}

\newtheorem{ecor}[equation]{Corollary}
\newtheorem{df}[equation]{D\'efinition}
\newtheorem{question}[equation]{Question}
\newenvironment{exemple}{
\vspace{2mm}
\refstepcounter{equation}\noindent\textbf{Exemple \theequation\ }}%
{\hfill\carrenoir\nolinebreak\vspace{2mm}}
\newenvironment{remarque}{
\refstepcounter{equation}\noindent\textbf{Remarque \thetheo\ }}
{\nolinebreak\vspace*{1mm}}

\makeatletter
\renewcommand\theequation{\thesection.\arabic{equation}}
\@addtoreset{equation}{section}
\makeatother
\newcommand{\carrenoir}{\rule{0.5em}{0.5em}}
\newenvironment{demo}[1][]{\textbf{D\'emonstration%
\ifthenelse{\equal{#1}{}}{}{ #1} :\ifthenelse{\equal{#1}{}}{ }{}}}
{\hfill\carrenoir\nolinebreak\vspace{2mm}}

\newcommand{\smat}[1]{\left(\begin{smallmatrix}#1\end{smallmatrix}\right)}

\newcommand{\oper}[2]{\newcommand{#1}{\mathop{\mathrm{#2}}\nolimits} }
\oper{\Ker}{Ker}
\oper{\Ima}{Im}
\oper{\Id}{Id}
\oper{\Vol}{Vol}
\oper{\diam}{diam}
\oper{\injrad}{inj}
\oper{\dimension}{dim}
\oper{\SL}{SL}
\oper{\GL}{GL}

\newcommand{\de}{\mathrm{ d }}
\newcommand{\trans}[1]{\vphantom{#1}^t#1}

\newcommand{\N}{\mathbb N}
\newcommand{\Z}{\mathbb Z}
\newcommand{\R}{\mathbb R}
\newcommand{\C}{\mathbb C}

\newcommand{\Q}{\mathbb Q}
\DeclareFontFamily{OML}{eur}{\skewchar\font127}
\DeclareFontShape{OML}{eur}{m}{n}{<5> <6> <7> <8> <9> gen * eurm
<10><10.95><12><14.4><17.28><20.74><24.88>eurm10}{}
\DeclareSymbolFont{greek}{OML}{eur}{m}{n}
\DeclareMathSymbol{\codiff}{\mathord}{greek}{"0E}
\DeclareMathSymbol{\prodint}{\mathord}{greek}{"13}

\title{Effondrement, spectre et propriétés diophantiennes des flots
riemanniens}
\author{Pierre Jammes}
\date{}
\begin{document}
\maketitle
{\small 
\textsc{Résumé.---}
On étudie le comportement des premières valeurs propres du laplacien 
agissant sur les formes différentielles lors d'un effondrement adiabatique 
d'un flot riemannien $\mathcal F$ sur une variété compacte $M$. Le nombre 
de petites valeurs
propres peut alors se calculer en fonction de la cohomologie basique
de $\mathcal F$, et on donne des critères spectraux pour l'annulation
des classes d'\'Alvarez et d'Euler du flot. En outre, on définit un
invariant de nature diophantienne du flot qui est lié au comportement
asymptotique des petites valeurs propres. Un appendice est consacré aux
propriétés arithmétiques des flots riemanniens.

Mots-clefs : effondrements, formes différentielles, laplacien, 
petites valeurs propres, flots riemanniens, approximations diophantiennes.

\medskip
\textsc{Abstract.---}
We study the behavior of the first eigenvalues of the Hodge Laplacian acting
on differential forms under adiabatic collapsing of a riemannian flow 
$\mathcal F$ on a closed manifold $M$. We show that the number of small
eigenvalues is related to the basic cohomology of $\mathcal F$, and
give spectral criteria for the vanishing of the \'Alvarez class
and the Euler class of the flow. We also define a diophantine
invariant of the flow which is related to the asymptotical behavior
of the small eigenvalues. An appendix is devoted to arithmetic
properties of riemannian flows.


}

\selectlanguage{english}
\renewcommand{\theequation}{\arabic{equation}}

\section*{Summary}
The aim of this article is the study the behavior of the first eigenvalues
of the Hodge Laplacian on certain collapsing manifolds. Most works on
his subject deal with collapsing of a fibre bundles on their basis. In this
article, we consider a diferent topological structure on the manifold,
namely a foliation. More precisely, we focus on the case of a riemannian
flow, that is a 1-dimensional foliation equiped with a transversal
riemannian structure: there exist a metric $g$ such that if we write
$g=g_H\oplus g_V$, where $g_V$ is the component of the metric tangential
to the leaves and $g_H$ the orthogonal component, $g_H$ is constant along
the flow. Such a metric is called ``bundle-like''. For a bundle-like metric 
$g$, we define the adiabatic collapsing associated
to $g$ by $g_\varepsilon=\varepsilon^2g_H\oplus g_V$. The curvature of
$g_\varepsilon$ is known to be bounded as $\varepsilon$ tends to $0$.

 The behavior of the firsts eigenvalues of the Hodge laplacian is
related to the basic cohomology $H^*(M/\mathcal F)$ of the flow, which
is the cohomology of the space of the basic forms $\Omega^*(M/\mathcal F)
=\{\omega\in\Omega^*(M),\ \prodint_X\omega=\prodint_X\de\omega=0\}$, where
$X$ is a vector field tangent to the flow. This
cohomology is finite dimensional. If $0<\lambda_{p,1}(M,g)\leq 
\lambda_{p,2}(M,g)\leq \ldots$ denotes the positive spectrum of the laplacian
on $p$-forms, we prove the following:

\begin{etheo}\label{sum:th1}
Let $M$ a compact $n$-manifold, $\mathcal F$ a riemannian flow on $M$, 
and $(g_\varepsilon)$ the adiabatic collapsing of a bundle-like metric $g$.
 Then :
\begin{enumerate}
\item The number of positive small eigenvalues on $p$-forms
for the collapsing $(g_\varepsilon)$ is
$$m_p=\dimension H^p(M/\mathcal F)+\dimension H^{n-p}(M/\mathcal F)-b_p(M);$$
\item there is a constant $c(g,\mathcal F)>1$ such that for all $p$, $k$,
and $\varepsilon$ such that
$1\leq k\leq m_p$ and $0<\varepsilon<1$, 
$c^{-1}\cdot\varepsilon^2\leq
\lambda_{p,k}(M,g_\varepsilon)\leq c\cdot\varepsilon^2$.
\end{enumerate}
\end{etheo}

The existence of small eigenvalues is related to two cohomological
invariants of the flow, namely the \'Alvarez class and the Euler class.
According to \cite{do98}, there is a bundle-like metric such that
if $X$ is a unitary vector field tangent to the flow $\mathcal F$ and 
$\chi=X^\flat$ its characteristic form, the basic component of the 1-form
$\kappa=\mathcal L_X\chi$ is closed. Its basic cohomology class is denoted
$[\kappa]\in H^1(M/\mathcal F)$, and is an invariant of the flow 
(see~\cite{al92}). It is an obstruction for the flow to be isometric.
The Euler form of $\mathcal F$ is defined by $e=\de\chi-\kappa\wedge\chi$
and induces a class $[e]\in H^2_{-\kappa}(M/\mathcal F)$, where
$H^*{-\kappa}(M/\mathcal F)$ is the twisted basic cohomology defined
by the twisted differential operator $\de_{-\kappa}=\de+\kappa\wedge$ 
(see~\cite{rp01a}). We can deduce from the theorem~\ref{sum:th1} the 
following vanishing criteria for the \'Alvarez and Euler class:
\begin{ecor}
If there is a degree $p$ such that $\lambda_{p,1}(M,g_\varepsilon)$ 
tends to 0 when
$\varepsilon\to0$, then $[e]\neq0$.

The eigenvalue $\lambda_{1,1}(M,g_\varepsilon)$ tends to 0 as
$\varepsilon\to0$ if and only if $[\kappa]=0$ and $[e]\neq0$.
\end{ecor}

The theorem~\ref{sum:th1} gives bounds of the small eigenvalues depending
on the parameter $\varepsilon$, wich is proportional to the volume of the
manifold. In the second part of the article, we study the behavior of
small eigenvalues relatively to the Gromov-Hausdorff distance between
the manifold $(M,g_\varepsilon)$ and the limit space of the collapsing.
The closure of a leaf is known to be a torus on which the flow is linear
(cf.~\cite{ca84a}, \cite{ca84b}). We prove that the dimension of this torus
and this linear flow are the same for almost all closures of leaves
(lemma~\ref{geom:lem0}). Hence, the slope of the linear flow on
closures of generic leaves is an invariant of the flow. If this
closure is two dimensional, the slope of the flow is an irrational number
$\alpha\in\R\backslash\Q$, and we can define an irrationnality exponent
for the flow by 
$\mu(\mathcal F)=\sup\{\nu,\ |\alpha-\frac pq|<\frac1{q^\nu} \textrm{ has
infinitely many solutions}\}$.

\begin{etheo}
Let $\mathcal F$ be a riemannian flow on a compact manifold $M$
such that the closure of generic leaves of $\mathcal F$ is of dimension~2.
Then, for all adiabatic collapsing $(g_\varepsilon)$ of $\mathcal F$, 
$$\liminf_{\varepsilon\to0}\frac{\ln\delta(M,\mathcal F,g_\varepsilon)}
{\ln\varepsilon}=\frac1{\mu(\mathcal F)},$$
where $\mu(\mathcal F)$ denotes the irrationnality exponent of $\mathcal F$,
and $\delta(M,\mathcal F,g_\varepsilon)$ the Gromov-Hausdorff distance 
between $(M,g_\varepsilon)$ and the limit space of the collapsing.
\end{etheo}
It follows that, under these hypothesis, we have
$\limsup_{\varepsilon\to0}\frac{\ln\lambda_{p,i}(M,g_\varepsilon)}
{\ln\delta(M,\mathcal F,g_\varepsilon)}=2\mu(\mathcal F)$.

When the dimension $k$ of the closure of generic leaves is greater than~3,
the slope of the flow is a vector $\alpha\in\R^{k-1}$. If $\alpha$
is badly approximable, that is
$\|\alpha-\frac pq\|\geq\frac c{q^{1+\frac1{k-1}}}$
for some $c(\alpha)>0$ and for all  $q\in\Z^*$, $p\in\Z^{k-1}$ we have
the following :

\begin{etheo}
Let $M$ a compact manifold, and $\mathcal F$ a riemannian flow on $M$
If the generic closure of leaves of $\mathcal F$ are $k$ dimensional
and the slope of the flow is badly approximable, then
for all adiabatic collapsing $(g_\varepsilon)$ of $\mathcal F$,
there is a constant $c>0$ such that
\begin{equation}
\frac1c\cdot\varepsilon^{\frac1k}\leq\delta(M,\mathcal F,g_\varepsilon)
\leq c\cdot\varepsilon^{\frac1k}
\end{equation}
for all $\varepsilon\in]0,1]$.
\end{etheo}
As a corollary, we have 
$\frac1c\cdot\delta(M,\mathcal F,g_\varepsilon)^{2k}\leq
\lambda_{p,i}(M,g_\varepsilon)\leq
c\cdot\delta(M,\mathcal F,g_\varepsilon)^{2k}$.

 The last part of the article is devoted to examples of collapsing
riemannian flows. In particular, we exhibit flows with prescibed 
$\mu(F)$ and with badly approximable slope.
\selectlanguage{french}
\renewcommand{\theequation}{\thesection.\arabic{equation}}

\section{Introduction}
On sait qu'à diamètre borné et courbure de Ricci minorée, la première
valeur propre du laplacien agissant sur les fonctions d'une variété
compacte est uniformément minorée. Dans \cite{cc90}, B.~Colbois
et G.~Courtois ont montré que ce résultat ne s'étendait pas au
laplacien agissant sur les formes différentielles et que même avec
une hypothèse de courbure sectionnelle bornée, on pouvait trouver
des variétés admettant une suite de métriques telle que la première
valeur propre du laplacien tende vers 0. Ils montrent en outre
qu'à courbure et diamètre bornés, si une valeur propre tend vers zéro, 
alors le volume ---~ou de manière équivalente 
le rayon d'injectivité~--- de la variété tend vers 0, c'est-à-dire
qu'elle s'effondre. Ces résultats motivent le problème suivant :
\begin{question}\label{intro:q}
À quelles conditions une variété qui s'effondre admet-elle une
ou plusieurs petites valeurs propres~? Peut-on estimer à quelle vitesse
ces valeurs propres tendent vers zéro par rapport au volume ou au rayon
d'injectivité~?
\end{question}
Ces questions ont déjà fait l'objet de plusieurs travaux (\cite{cc00},
\cite{lo02}, \cite{ja03}, \cite{ja04}, voir \cite{ja05} pour une présentation
synthétique de ces résultats), mais les situations étudiées
sont celles d'une variété $M$ tendant pour la distance de
Gromov-Hausdorff vers une variété lisse $N$ de dimension inférieure, $M$
ayant alors une structure de fibré sur $N$, alors que l'espace limite
d'une variété qui s'effondre est en général une variété stratifiée.
Les seuls exemples connus de petites valeurs propres dans 
le cas où la variété $M$ tend vers un espace métrique qui n'est pas une 
variété lisse sont des exemples de fibrés de Seifert (\cite{cc90}),
situation très proche des fibrés en cercles étudiés en détail dans \cite{cc00}.

On va s'intéresser ici à des situations plus générales que les fibrés,
à savoir les feuilletages, en se restreignant au cas plus simple
où les feuilles sont de dimension~1, c'est-à-dire aux flots, et on
considérera des effondrements obtenus en faisant varier la métrique le long 
des feuilles. Le choix des feuilletages est en particulier motivé par le 
fait qu'il permet d'exhiber facilement des exemples de petites valeurs propres 
dans le cas où la variété s'effondre sur un espace singulier (voir
ci-dessous), et il fera en outre apparaître des différences notables
avec le cas des espaces limites lisses (voir remarque~\ref{intro:rq2}). 
La restriction aux flots s'explique par le fait que même 
dans le cas des fibrés, la seule situation vraiment élucidée est celle 
où la fibre est un cercle. Le fait que la variété s'effondre à courbure 
bornée impose aussi certaines contraintes sur le feuilletage que nous 
détaillerons plus loin.

Il s'avère que l'exemple le plus simple d'effondrement sur un espace
singulier, qui est l'effondrement d'un flot isométrique, fournit 
un exemple de petite valeur propre. Et contrairement aux exemples donnés
dans \cite{cc90}, on peut faire en sorte que les feuilles ne soient pas
compactes:
\begin{exemple}\label{intro:ex}
On considère sur la sphère $S^3=\{(a,b)\in\C^2,\ |a|^2+|b|^2=1\}$ l'action
isométrique du tore $T^2$ définie par $(\theta_1,\theta_2)\cdot
(a,b)=(e^{i\theta_1}a,e^{i\theta_2}b)$. Si on se donne un irrationnel
$\alpha\in\R\backslash\Q$, on peut lui associer un flot sur $S^3$ par 
le plongement de $\R$ dans $T^2$ défini par $t\mapsto(t,\alpha t)$ en 
considérant l'action induite sur $S^3$.

Pour construire l'effondrement, on décompose la métrique  canonique $g$
en la somme $g=g_H\oplus g_V$ d'une composante verticale, tangente au flot, 
et d'une composante horizontale $g_H$ sur l'espace orthogonal au flot.
On définit alors la famille de métriques $g_\varepsilon=g_H\oplus 
\varepsilon^2g_V$. Comme $\alpha$ est irrationnel, l'adhérence de chaque
feuille est une orbite de l'action de $T^2$, et donc $(S^3,g_\varepsilon)$
tend vers $S^3/T^2=[0,\frac\pi2]$ pour la distance de Gromov-Hausdorff
quand $\varepsilon$ tend vers 0.

On considère la 1-forme différentielle $\omega=\frac X{|X|^2}$, où $X$ est
le champ de vecteur associé au flot. Sa codifférentielle $\codiff\omega$
est nulle : en effet, $\codiff\omega$ est une fonction invariante par le
flot, donc $\de\codiff\omega$ est partout orthogonale à $\omega$,
et donc $\|\codiff\omega\|^2=(\omega,\de\codiff\omega)=0$. En outre,
sa différentielle vérifie $i_X\de\omega=\mathcal L_X\omega-\de i_X\omega=0$.
Le quotient de Rayleigh de $\omega$ d'écrit donc
$R(\omega)=\frac{\|\de\omega\|^2}{\|\omega\|^2}$, et tend vers 0 quand
$\varepsilon$ tend vers zéro, car $\|\omega\|^2$ tend vers l'infini et
$\|\de\omega\|^2$ reste constant. On est alors assuré
qu'une valeur propre du laplacien agissant sur $\Omega^1(S^3)$ tend vers
zéro quand on effondre la sphère, et cette valeur propre est non nulle
puisque $b_1(S^3)=0$.
\end{exemple}

On va chercher dans cet article à comprendre dans quelle mesure
cet exemple se généralise aux autres flots. On considère donc un flot
$\mathcal F$ ---~c'est-à-dire un feuilletage orientable de dimension 1,
\emph{a priori} sans paramètre~---
sur une variété $M$ et si on se donne une métrique $g$ sur $M$, on définit
comme dans l'exemple \ref{intro:ex} une famille de métriques
$(g_\varepsilon)$ en décomposant $g$ sous la forme $g_H\oplus g_V$, où
$g_V$ est la métrique le long des feuilles et $g_H$ la composante
de la métrique $g$ orthogonale aux feuilles, et en posant
$g_\varepsilon=g_H\oplus\varepsilon^2g_V$ pour $\varepsilon\leq1$. 
On appellera «~effondrement
adiabatique associé à $g$~» la famille de métrique ainsi construite.

Il est rare qu'une telle déformation de la métrique maintienne 
la courbure bornée, même avec l'hypothèse que les feuilles sont
de dimension~1. On supposera donc en outre que le flot $\mathcal F$ est 
riemannien, c'est-à-dire qu'il admet une métrique, dite «~quasi-fibrée~», 
dont la composante $g_V$ est invariante par le flot (cette 
propriété ne dépend pas du paramétrage choisi sur le flot). En effet,
Y.~Carrière a montré (\cite{ca84b}) que si $g$ est une telle métrique, 
la courbure sectionnelle de $g_\varepsilon$ est uniformément bornée par 
rapport à $\varepsilon$, cette propriété étant fausse pour les
feuilletages riemanniens de dimension plus grande (même en supposant par
exemple que le feuilletage est défini par une action localement libre
de $\R^k$). Les flots riemanniens 
fournissent ainsi un grand nombre d'exemples d'effondrements sur des 
espaces singuliers (nous en rappelerons certains dans la section~\ref{ex}). 
Rappelons aussi une autre propriété de ces métriques, 
que nous n'utiliserons pas ici mais qui intervient 
dans d'autres contextes (voir par exemple \cite{gh83}, \cite{mo05}, 
\cite{ma08}
et les références qui y sont données) et qui motive l'étude des flots 
riemanniens en général : si on se donne un champ d'hyperplans
sur une variété riemannienne, il est totalement géodésique
si et seulement si la métrique est quasifibrée pour le flot orthogonal,
ce flot étant alors riemannien. 

Sous ces hypothèses, nous allons mettre en évidence quels sont les
points communs et surtout les différences avec les effondrements
de fibrés. En particulier, nous tacherons d'éclairer les liens entre 
comportement du spectre, géométrie de l'effondrement et dynamique du flot 
qui n'étaient qu'esquissés dans \cite{ja04}. 

Si $X$ est un champ de vecteur tangent au flot, on définit la cohomologie 
basique $H^*(M/\mathcal F)$ du flot comme étant la cohomologie de l'espace 
des formes basiques $\Omega^*(M/\mathcal F)=\{\omega\in\Omega^*(M),\ 
\prodint_X\omega=\prodint_X\de\omega=0\}$. On sait (\cite{eksh85}) que si 
le flot est riemannien, cette cohomologie
est de dimension finie. On va montrer qu'on peut calculer le nombre
de petites valeurs propres pour l'effondrement $(g_\varepsilon)$
en fonction de cette cohomologie. Pour tout 
$p$, on notera $\lambda_{p,0}(M,g)$ la valeur propre nulle pour les 
$p$-formes, si elle existe, et $\lambda_{p,1}(M,g)\leq \lambda_{p,2}(M,g)\leq
\ldots$ les valeurs propres non nulles, en les répétant s'il
y a multiplicité.
\begin{theo}\label{intro:th1}
Soit $(M,g)$ une variété riemannienne compacte de dimension $n$, 
$\mathcal F$ un flot
riemannien sur $M$ et $(g_\varepsilon)$ l'effondrement adiabatique 
associé à une métrique $g$ quasi-fibrée pour $\mathcal F$. Alors :
\begin{enumerate}
\item le nombre de petites valeurs propres non nulles sur les $p$-formes 
pour l'effondrement $(g_\varepsilon)$ est 
$$m_p=\dimension H^p(M/\mathcal F)+\dimension H^{n-p}(M/\mathcal F)-b_p(M)\ ;$$
\item il existe une constante $c(g,\mathcal F)>1$ telle que pour tout $p$,
tout $1\leq k\leq m_p$ et tout $0<\varepsilon<1$, on a 
$c^{-1}\cdot\varepsilon^2\leq
\lambda_{p,k}(M,g_\varepsilon)\leq c\cdot\varepsilon^2$.
\end{enumerate}
\end{theo}
\begin{remarque} Le théorème \ref{intro:th1} s'applique en particulier
aux fibrés en cercles, ce qui permet de retrouver une partie des résultats de 
\cite{cc00}.
\end{remarque}

\begin{remarque}\label{intro:rq2}
 En revanche, contrairement à ce qui se passe pour
les effondrements de fibrés sur leur base, ce n'est pas la cohomologie
de l'espace limite qui intervient ici mais celle de la structure 
transverse. Cela contraste aussi avec les exemples d'effondrements à 
courbure minorée étudiés par J.~Lott dans \cite{lo04} pour lesquels
le nombre de petite valeurs propres se calcule à l'aide de la 
cohomologie de l'espace limite. La structure
de fibré sur l'espace limite ne semble donc pas pertinente pour étudier
dans toute sa généralité le problème posé par la question~\ref{intro:q}.
\end{remarque}

On peut déduire du théorème \ref{intro:th1} un lien entre les propriétés 
spectrales du flot et
deux invariants cohomologiques, les classes d'\'Alvarez et d'Euler, dont
nous allons rappeler les définitions :

On se donne une métrique $g$ quasi-fibrée pour le flot, un champ
de vecteur $X$ unitaire et tangent au flot. La forme caractéristique du
flot est $\chi=X^\flat$ et sa forme de courbure moyenne est $\kappa
=\mathcal L_X\chi$. Dans \cite{al92}, J.-A.~\'Alvarez~L\'opez a montré que
la composante basique de cette forme de courbure est fermée, et que
la classe de cohomologie basique de cette composante ne dépend pas de
la métrique. On appelle \emph{classe d'\'Alvarez} cette classe de
cohomologie et on la note $[\kappa]\in H^1(M/\mathcal F)$. On sait
que cette classe est nulle si et seulement si le flot est géodésible, ou
de manière équivalente qu'il est isométrique ---~c'est-à-dire qu'il 
existe une métrique sur $M$ et un paramétrage du flot 
tels que le flot agisse par isométrie.

On définit la forme d'Euler du flot par $e=\de\chi-\kappa\wedge\chi$.
Cette forme est fermée pour l'opérateur différentiel tordu $\de_{-\kappa}$
défini par $\de_{-\kappa}\omega=\de\omega+\kappa\wedge\omega$,
et il existe une métrique quasi-fibrée pour laquelle cette forme est basique.
La forme $e$ représente donc une classe de cohomologie $[e]$ dans la 
cohomologie basique tordue $H^2_{-\kappa}(M/\mathcal F)$ définie comme 
la cohomologie de $\Omega^*(M/\mathcal F)$ pour la différentielle tordue
$\de_{-\kappa}$. La classe $[e]$ est appelée \emph{classe d'Euler} du flot 
$\mathcal F$.
Elle dépend de la métrique, mais pas le fait qu'elle soit nulle:
on peut montrer (\cite{rp01a}, \cite{rp01b}) que $[e]=0$ si et seulement 
s'il existe un feuilletage transverse à $\mathcal F$ dont la torsion est 
basique.

On peut alors écrire :
\begin{ecor}\label{intro:cor1}
S'il existe $p$ tel que $\lambda_{p,1}(M,g_\varepsilon)$ tend vers 0 quand
$\varepsilon$ tend vers 0, alors $[e]\neq0$.

La valeur propre $\lambda_{1,1}(M,g_\varepsilon)$ tend vers 0 quand 
$\varepsilon$ tend vers 0 si et seulement si $[\kappa]=0$ et $[e]\neq0$.
\end{ecor}
\begin{remarque}
Si $H^1(M)=\{0\}$, par exemple si $M$ est simplement connexe, 
on a nécessairement $[\kappa]=0$ et $[e]\neq0$ (voir
remarque \ref{coh:rq2} dans la section \ref{coh}). Cette situation 
est illustrée par l'exemple \ref{intro:ex}.
\end{remarque}

 Dans le cas où le flot est isométrique, on a donc une condition
nécessaire et suffisante sur le spectre ---~à savoir l'absence de petites
valeurs propres~--- pour que la classe d'Euler 
s'annule, comparable à celle qui est déjà connue dans le cas des
fibrés en cercle (\cite{cc00}). Cette condition ne se généralise cependant
pas au cas des flots non isométriques : on donnera au paragraphe
\ref{ex:euler} un exemple de flot dont la classe d'Euler est non nulle
mais dont l'effondrement ne produit pas de petite valeur propre.
On verra aussi dans la section \ref{coh} que ce phénomène peut s'exprimer 
en termes cohomologiques (remarque \ref{coh:nul}).

Le théorème \ref{intro:th1} donne l'ordre de grandeur de la vitesse
à laquelle les petites valeurs propres tendent vers zéro par rapport 
à $\varepsilon$, c'est-à-dire par rapport au volume. Dans le cas des
fibrés en cercles, le volume, le rayon d'injectivité et la distance de
Gromov-Hausdorff sont tous du même ordre au cours de l'effondrement; 
ce n'est pas le cas pour des effondrements plus généraux, et en
particulier pour les effondrements de flots riemanniens dont les feuilles
ne sont pas toutes compactes. On a construit dans \cite{ja04} des exemples 
d'effondrement de flots isométriques pour lesquels le comportement 
asymptotique de la première valeur propre non nulle par rapport au rayon 
d'injectivité dépend d'une propriété de nature dynamique du flot et qui 
s'exprime en termes d'approximations diophantiennes. 
Ces exemples ne se généralisent pas à 
tous les flots riemanniens (voir la remarque \ref{geom:rq}), 
mais on peut définir un invariant diophantien lié au comportement
asymptotique des petites valeurs propres par rapport à la distance de
Gromov-Hausdorff entre $M$ et l'espace limite de l'effondrement,
dans le cas d'effondrements adiabatiques.

On sait que les adhérences des feuilles d'un flot riemannien sont des tores
et qu'en restriction à l'adhérence d'une feuille, le flot est conjugué
à un flot linéaire (cf. \cite{ca84a}). On peut donc décrire le flot
sur une adhérence de feuille difféomorphe à $T^k$ par un
vecteur $(1,\alpha)$ avec $\alpha\in\R^{k-1}$. Ce vecteur n'est
défini qu'à un difféomorphisme linéraire du tore près, mais cela 
n'influe pas sur ses propriétés diophantiennes. On peut en fait montrer qu'il
existe un entier $k$ et un flot linéaire $\mathcal F_0$ sur $T^k$ tels que
sur presque toutes les adhérences des feuilles du flot $\mathcal F$,
la restriction de $\mathcal F$ est conjuguée à $\mathcal F_0$
(voir lemme~\ref{geom:lem0} pour un énoncé précis), ce qui permet 
de considérer les propriété diophantiennes de $\mathcal F_0$ comme
des propriétés globales du flot $\mathcal F$. On parlera 
d'\emph{adhérence générique} pour désigner les éléments de la classe 
de conjugaison de $\mathcal F_0$, et de \emph{feuille générique}
pour désigner les feuilles dont l'adhérence est générique. On peut aussi 
montrer que la dimension maximale des adhérences des feuilles est
celle de l'adhérence générique.

 Supposons dans un premier temps que l'adhérence générique du flot 
$\mathcal F$ est de dimension~2. En restriction à une adhérence
générique le flot s'identifie à un flot linéraire
sur $\R^2/\Z^2$ induit par un vecteur qui peut s'écrire $(1,\alpha)$ dans
une base de $\Z^2$, avec 
$\alpha\in\R\backslash\Q$. 
L'exposant d'irrationalité (appelé aussi \emph{mesure d'irrationalité}) 
de $\alpha$ est défini par 
\begin{equation}
\mu(\alpha)=\sup\left\{\nu,\ |\alpha-\frac pq|<\frac1{q^\nu} 
\textrm{ a une infinité de solutions }(p,q)\in\Z^2\right\}.
\end{equation}
Rappelons que l'exposant d'irrationalité est au moins égal à deux 
sur les irrationnels, qu'il vaut exactement 2 pour presque tout réel,
en particulier pour les réels algébriques, et que les réels pour lesquels
il est infini sont exactement les nombres de Liouville.

 Cet exposant ainsi associé à un flot linéaire sur $T^2$ ne dépend pas de 
la base de $\Z^2$ choisie
et, en vertu du lemme \ref{geom:lem0}, on peut définir un invariant global
du flot en considérant l'exposant d'irrationalité du flot sur une adhérence 
générique :
\begin{df} Soit $\mathcal F$ un flot riemannien sur une variété
compacte $M$ dont les adhérences génériques sont de dimension 2. 
On appelle exposant d'irrationalité du flot, noté $\mu(\mathcal F)$, l'
exposant d'irrationalité du vecteur directeur du flot restreint à une
adhérence générique.
\end{df}

On a alors, en notant $\delta(M,\mathcal F,g_\varepsilon)$ la distance
de Gromov-Hausdorff entre $(M,g_\varepsilon)$ et l'espace limite de
l'effondrement :
\begin{theo}\label{intro:th2}
Soit $M$ une variété riemannienne compacte et $\mathcal F$ un flot riemannien 
sur $M$ dont la dimension maximale des adhérence des feuilles est 2. 
Alors, pour tout effondrement adiabatique $(g_\varepsilon)$ du flot, on a
\begin{equation}
\liminf_{\varepsilon\to0}\frac{\ln\delta(M,\mathcal F,g_\varepsilon)}
{\ln\varepsilon}=\frac1{\mu(\mathcal F)},
\end{equation}
où $\mu(\mathcal F)$ est l'exposant d'irrationalité du flot.
\end{theo}
\begin{remarque} On en déduit que sous les hypothèses du théorème 
\ref{intro:th2}, on a
\begin{equation}
\limsup_{\varepsilon\to0}\frac{\ln\lambda_{p,i}(M,g_\varepsilon)}
{\ln\delta(M,\mathcal F,g_\varepsilon)}=2\mu(\mathcal F),
\end{equation}
pour tout entier $p$ et tout $0<i\leq m_p$.
\end{remarque}

\begin{remarque}On verra dans la section \ref{ex} (paragraphe \ref{ex:ex})
que pour tout
$\mu\in[2,+\infty]$, on peut construire des flots dont l'exposant
d'irrationalité est $\mu$.
\end{remarque}

La démonstration du théorème \ref{intro:th2} ne se généralise pas 
quand la dimension des adhérences des feuilles dépasse 2. Cependant, si la
direction du flot sur une adhérence générique
est difficilement approchable, c'est-à-dire induite par un vecteur
$(1,\alpha)\in\R^k$ tel qu'il existe $c(\alpha)$>0 pour lequel
\begin{equation}
\|\alpha-\frac pq\|\geq\frac c{q^{1+\frac1{k-1}}}
\end{equation}
pour tout $q\in\Z^*$ et $p\in\Z^{k-1}$,
on peut déterminer précisément le comportement
asymptotique de $\delta(M,\mathcal F,g_\varepsilon)$ :
\begin{theo}\label{intro:th3}
Soit $M$ une variété riemannienne compacte et $\mathcal F$ un flot riemannien
sur $M$. Si les adhérences génériques du flot sont de dimension $k$ et
que sur ces adhérences, la
direction du flot est difficilement approchable, alors pour
tout effondrement adiabatique $(g_\varepsilon)$ du flot, il existe 
une constante $c>0$ telle que
\begin{equation}
\frac1c\cdot\varepsilon^{\frac1k}\leq\delta(M,\mathcal F,g_\varepsilon)
\leq c\cdot\varepsilon^{\frac1k}
\end{equation}
pour tout $\varepsilon\in]0,1]$.
\end{theo}
\begin{remarque} On peut en déduire un résultat sur le spectre comme pour 
le théorème \ref{intro:th2}: sous les hypothèses du théorème \ref{intro:th3},
il existe une constante $c>0$ telle que
\begin{equation}
\frac1c\cdot\delta(M,\mathcal F,g_\varepsilon)^{2k}\leq
\lambda_{p,i}(M,g_\varepsilon)\leq
c\cdot\delta(M,\mathcal F,g_\varepsilon)^{2k},
\end{equation}
pour tout entier $p$, tout $0<i\leq m_p$ et tout $\varepsilon\in]0,1]$. 
Notons aussi que dans
\cite{ja04}, on ne donnait des exemples de valeurs propres de 
l'ordre de $\delta(M,\mathcal F,g_\varepsilon)^{2k}$ uniquement en
degré inférieur à 3. Ici, on peut effectivement avoir de tels
valeurs propres en tout degré (voir paragraphes \ref{ex:spheres} et
\ref{ex:ex}).
\end{remarque}

On verra dans la section \ref{ex} comment construire des flots dont
la direction est difficilement approchable (paragraphes \ref{ex:ex} et
\ref{ex:diff}). On peut aussi énoncer le
critère suivant (voir remarque \ref{ap:rq1} de l'appendice):
\begin{proposition}\label{intro:pr}
Si $\mathcal F$ est un flot riemannien non isométrique dont
la dimension de l'adhérence générique est un nombre premier, alors 
la direction du flot est difficilement approchable.
\end{proposition}
\begin{remarque}
Dans le cas ou l'adhérence générique est de dimension 2, le fait que
la direction du flot $\mathcal F$ soit difficilement approchable implique 
que $\mu(\mathcal F)=2$. L'exposant d'irrationalité ne peut donc être 
prescrit que dans le cadre des flots isométriques. On a là une illustration 
du caractère rigide des flots riemanniens non isométriques.
\end{remarque}

Le théorème~\ref{intro:th1} et le corollaire~\ref{intro:cor1} 
seront démontrés dans la section~\ref{coh}. La section~\ref{geom} 
sera consacrée à l'étude de $\delta(M,\mathcal F,g_\varepsilon)$ 
et à la démonstration des théorèmes~\ref{intro:th2} et~\ref{intro:th3}.
Dans la section~\ref{ex}, nous présenterons plusieurs exemples. Nous
terminerons par un appendice sur les propriétés arithmétiques des
flots riemanniens qui nous serons utiles pour construire certains exemples
et pour établir le proposition~\ref{intro:pr}.

\section{Cohomologie et spectre des flots riemanniens}\label{coh}
\subsection{Suite spectrale des flots riemanniens}
Les démonstrations du théorème \ref{intro:th1} et du corollaire
\ref{intro:cor1} s'appuient essentiellement sur les
résultats de J.~A.~\'Alvarez~L\'opez et Y.~A.~Kordyukov qui généralisent
dans \cite{alk00} des résultats antérieurs de Mazzeo, Melrose et Forman
(\cite{mm90}, \cite{fo95}) en montrant qu'on peut calculer le nombre de 
petites valeurs propres lors d'un effondrement adiabatique (sans 
hypothèse sur la courbure)
d'un feuilletage riemannien à l'aide de la suite spectrale différentielle
du feuilletage. Cette suite spectrale apparaît pour la première
fois dans les travaux de K.~S.~Sarkaria \cite{sa74} \cite{sa78}.
Nous renvoyons au chapitre 4 de \cite{to97} ou à \cite{alk00} pour une 
définition de la suite spectrale différentielle 
d'un feuilletage riemannien. Afin de ne pas alourdir le texte, nous 
noterons $E_1^{*,*}$ ce qui est noté $\hat{E}_1^{*,*}$ dans \cite{alk00}:
\begin{theo}[\cite{alk00}]\label{coh:alk1}
Soit $\mathcal F$ un feuilletage riemannien sur une variété compacte
$M$, $g$ une métrique sur $M$ et $(g_\varepsilon)$ l'effondrement
adiabatique associé. On a alors, pour tous entiers $p$ et $k\geq0$, 
en notant $E^{*,*}_k$ la suite spectrale du feuilletage et en posant
$E^p_k=\displaystyle\bigoplus_{a+b=p} E^{a,b}_k$ :
$$\dimension E^p_{k+1}=\#\{i\geq1,\ \lambda_{p,i}(M,g_\varepsilon)\in
O(\varepsilon^{2k})\textrm{ \emph{quand} }\varepsilon\to0\}+b_p(M).$$
\end{theo}
Pour obtenir une minoration des petites valeurs propres comme celle
du théorème \ref{intro:th1}, on doit compléter ce résultat :
\begin{theo}\label{coh:alk2}
Sous les hypothèses du théorème \ref{coh:alk1}, il existe $c>0$ dépendant
de $g$ mais pas de $\varepsilon$ tel que pour tout entier $p\geq0$ et tout
entier $i>\dimension E^p_{k+1}-b_p(M)$, on a
$$c\cdot\varepsilon^{2k}\leq\lambda_{p,i}(M,g_\varepsilon)$$
pour tout $\varepsilon\leq1$.
\end{theo}
Nous utiliserons ce résultat pour une métrique $g$ quasi-fibrée, mais 
il est vrai pour une métrique quelconque, comme le théorème~\ref{coh:alk1}. 
Nous allons donc le montrer dans le cas général.

\begin{demo} Nous allons utiliser le fait (démontré dans \cite{alk00})
qu'il existe une suite décroissante de sous-espaces de $\Omega^*(M)$
\begin{equation}
\Omega^*(M)\supset\mathcal H^*_1\supset\mathcal H^*_2\supset\mathcal H^*_3
\supset\ldots\supset\mathcal H^*_\infty
\end{equation} 
et un opérateur symétrique $\Delta'_\varepsilon:\Omega^*(M)\to\Omega^*(M)$
dont le spectre est le même que celui de $\Delta(M,g_\varepsilon)$, et
vérifiant les deux propriétés suivantes:
\begin{itemize}
\item $\mathcal H^*_k\simeq E^*_k$;
\item Si $\omega_i$ est une suite d'éléments de $\Omega^p(M)$ telle que
$\|\omega_i\|=1$ et $(\Delta'_\varepsilon\omega_i,\omega_i)\in 
o(\varepsilon_i^{2(k-1)})$ où $\varepsilon_i$ est une suite tendant vers $0$,
alors il existe une sous-suite convergente de $\omega_i$  
dont la limite est dans $\mathcal H^p_k$ (\cite{alk00}, théorème B).
\end{itemize}
S'il existe une suite $\varepsilon_i$ tendant vers $0$ et telle que 
$\lambda_{p,m-b_p(M)}(M,g_\varepsilon)\cdot\varepsilon_i^{-2k}$ tende vers
zéro pour un $m>\dimension E^p_{k+1}$, alors on peut construire
une suite $(\mathcal B_i)_{i\in\N}$ de famille de vecteurs $\mathcal B_i=
(\omega_i^{(j)})_{1\leq j\leq m}$ telle que 
chaque famille $\mathcal B_i$ soit orthonormée et que pour tout $j$,
on a $(\Delta'_\varepsilon\omega^{(j)}_i,\omega^{(j)}_i)\in
o(\varepsilon_i^{2k})$ quand $i\to+\infty$. On peut donc 
extraire de $(\mathcal B_i)$ une sous-suite convergente, dont la 
limite sera orthonormée et contenue dans $\mathcal H^p_{k+1}$. Mais alors
on a $\dimension\mathcal H^p_{k+1}\geq m$, ce qui contredit 
$\dimension\mathcal H^p_{k+1}=\dimension E^p_{k+1}<m$.

Le théorème B de \cite{alk00}, utilise comme hypothèse que la métrique
$g$, et donc l'effondrement $(g_\varepsilon)$, est quasi-fibrée. Il reste
donc à montrer qu'on peut étendre le résultat qu'on vient d'obtenir
à des métriques quelconques. Soit $g$ une métrique quasi-fibrée, $g'$
une métrique quelconque, et $(g_\varepsilon)$ et $(g'_\varepsilon)$ les
effondrements adiabatiques associés. Comme la variété $M$ est compacte,
il existe des constantes $\tau>0$ et $\tau'>0$ telles que
$\tau' g\leq g'\leq\tau g$. Comme l'ont remarqué \'Alvarez et Kordyukov
dans \cite{alk00} (section 4), on a alors 
$\tau' g_\varepsilon\leq g'_\varepsilon\leq\tau g_\varepsilon$ pour tout 
$\varepsilon>0$. Un théorème de J.~Dodziuk (\cite{do82}) assure alors que 
\begin{equation}
\frac1{\tau'}\left(\frac{\tau'}\tau\right)^{\frac{3n}2}
\lambda_{p,k}(M,g_\varepsilon)\leq\lambda_{p,k}(M,g'_\varepsilon)\leq
\frac1\tau\left(\frac{\tau\vphantom{'}}{\tau'}\right)^{\frac{3n}2}
\lambda_{p,k}(M,g_\varepsilon)
\end{equation}
pour tout entiers $k\geq0$ et $p\in[0,n]$, ce qui montre bien que
le comportement asymptotique du spectre lors de l'effondrement 
ne dépend pas de la métrique initiale choisie.
\end{demo}

Dans le cas où les feuilles sont de dimension 1, la suite spectrale
du feuilletage est assez simple: le complexe bigradué
$\Omega^{*,*}(M/\mathcal F)$ sur lequel est construit la suite se réduit
aux deux lignes $\Omega^{*,0}(M/\mathcal F)$ et $\Omega^{*,1}(M/\mathcal F)$,
$\Omega^{*,q}(M/\mathcal F)$ étant nul pour $q\geq2$. Le terme $E^{*,0}_2$ 
s'identifie à la cohomologie basique $H^*(M/\mathcal F)$ comme pour
tous les feuilletages, et le terme $E^{*,1}_2$ est isomorphe à 
la cohomologie tordue $H^*_\kappa(M/\mathcal F)$ qui est la cohomologie
des formes basiques pour l'opérateur différentiel tordu $\de_\kappa:
\omega\to\de\omega-\kappa\wedge\omega$ (\cite{kt84}, diagramme 2.18).
Comme $E^{*,q}_2$ est nul pour $q\geq2$, la suite dégénère au rang 3,
et on peut déduire du théorème \ref{coh:alk2} que toutes les petites
valeurs propres sont de l'ordre de $\varepsilon^2$ et que leur nombre
se calcule en fonction de $H^*(M/\mathcal F)$, $H^*_\kappa(M/\mathcal F)$
et $b_p(M)$:
\begin{fait} Le nombre de petites valeurs propres non nulles sur les 
$p$-formes pour un effondrement $(g_\varepsilon)$ associé à une
métrique quasi-fibrée $g$ est 
$m_p=\dimension H^p(M/\mathcal F)+\dimension H_\kappa^{p-1}(M/\mathcal F)
-b_p(M)$, et il existe une constante $c(g,\mathcal F)>1$ telle que pour 
tout $p$ et tout $1\leq k\leq m_p$, on a $c^{-1}\cdot\varepsilon^2\leq
\lambda_{p,k}(M,g_\varepsilon)\leq c\cdot\varepsilon^2$.
\end{fait}
\begin{demo} La suite spectrale $E_k^{*,*}$ du feuilletage dégénère au
rang 3, donc $\dimension E_k^p=b_p$ pour tout $k\geq3$. Les seules 
petites valeurs propres sont donc celles correspondant à $k=1$ dans les
théorèmes \ref{coh:alk1} et \ref{coh:alk2}. Leur nombre est 
\begin{eqnarray}
m_p&=&\dimension E_2^p-b_p\nonumber\\
&=&\dimension E_2^{p,0}+\dimension E_2^{p-1,1}-b_p\nonumber\\
&=&\dimension H^p(M/\mathcal F)+\dimension H^{p-1}_\kappa(M/\mathcal F)-b_p,
\end{eqnarray}
et les théorèmes \ref{coh:alk1} et \ref{coh:alk2} assurent l'existence d'une
constante $c$ telle que 
\begin{equation}
c^{-1}\cdot\varepsilon^2\leq\lambda_{p,k}(M,g_\varepsilon)
\leq c\cdot\varepsilon^2
\end{equation}
 pour tout 
$1\leq k\leq m_p$.
\end{demo}

La cohomologie basique ne vérifie en général pas la dualité de Poincaré,
mais il existe une relation de dualité avec la cohomologie tordue:
\begin{fait}[\cite{kt83}]\label{coh:dual}
Si $\mathcal F$ est un flot riemannien sur une variété compacte $M$
de dimension $n$, alors on a $H_\kappa^i(M/\mathcal F)\simeq H^{n-1-i}
(M/\mathcal F)$.
\end{fait}
On en déduit que la constante $m_p$ peut s'écrire
\begin{equation}
m_p=\dimension H^p(M/\mathcal F)+\dimension H^{n-p}(M/\mathcal F)-b_p(M),
\end{equation}
Ce qui achève la démonstration du théorème \ref{intro:th1}.
\subsection{Suite de Gysin des flots riemanniens}
Pour démontrer le corollaire \ref{intro:cor1},
on va préciser les rôles des classes $[\kappa]\in H^1(M/\mathcal F)$ et
$[e]\in H^2_{-\kappa}(M/\mathcal F)$ dans la cohomologie du flot. On
a déjà dit que la classe d'\'Alvarez est nulle si et seulement si
le flot est isométrique. Cette propriété du flot possède d'autres
caractérisations:
\begin{theo}[\cite{ca84a}, \cite{kt83}, \cite{ms85}, \cite{al92}]%
\label{coh:geod}
Soit $\mathcal F$ est un flot riemannien sur une variété 
compacte $M$ de dimension $n$. Les propriétés suivantes sont équivalentes :
\begin{enumerate}
\item le flot $\mathcal F$ est isométrique;
\item le flot $\mathcal F$ est géodésible;
\item la classe d'\'Alvarez du flot est nulle;
\item $H^{n-1}(M/\mathcal F)\neq 0$;
\item la cohomologie basique du flot vérifie la dualité de Poincaré.
\end{enumerate}
\end{theo}

La classe d'Euler $[e]$ intervient dans une suite exacte longue de
cohomologie, construite par J.~I.~Royo~Prieto et qui généralise la suite
de Gysin des fibrés en cercles :
\begin{theo}[\cite{rp01a}, \cite{rp01b}]\label{coh:gysin}
 Si $\mathcal F$ est un flot riemannien sur une variété compacte $M$,
on a la suite exacte longue
\begin{equation}
\cdots\to H^i(M/\mathcal F)\to H^i(M)\to H_\kappa^{i-1}(M/\mathcal F)
\stackrel{\wedge e}{\to}H^{i+1}(M/\mathcal F)\to\cdots
\end{equation}
\end{theo}
\begin{remarque} Dans \cite{rp01a}, le théorème \ref{coh:gysin} est
démontré en utilisant un argument à la Mayer-Vietoris. On peut aussi
le déduire de la suite spectrale du fibré comme pour les fibrés en
cercles, en s'appuyant sur le fait que la suite spectrale vérifie
$E_2^{*,q}=0$ pour $q\geq2$ et en suivant la même démarche que dans
\cite{bt82} (chapitre 3, \S 14).
\end{remarque}

\begin{remarque}\label{coh:rq2}
Il est bien connu que si $M$ est simplement connexe, alors le flot
est nécessairement isométrique : c'est par exemple un cas particulier
d'un résultat d'É.~Ghys (\cite{gh84}) et on peut le retrouver en remarquant
que $[\kappa]\in H^1(M/\mathcal F)\subset H^1(M)=\{0\}$. La suite
de Gysin montre que l'hypothèse $H^1(M)=\{0\}$ implique aussi $[e]\neq0$.
En effet, on a alors la suite exacte 
\begin{equation}
\{0\}\to H_\kappa^0(M/\mathcal F) \stackrel{\wedge e}{\to}H^2(M/\mathcal F)\to
\end{equation}
et comme $[k]=0$, on a $H_\kappa^0(M/\mathcal F)=H^0(M/\mathcal F)\simeq\R$.
L'application $\wedge e$ étant injective, elle est donc aussi non triviale. 
Par conséquent $[e]\neq0$.
\end{remarque}

Le corollaire découle des résultats énoncés précédemment.\\ 
\begin{demo}[du corollaire \ref{intro:cor1}]
Si la classe d'Euler $[e]$ est nulle, la suite de Gysin se décompose
en suites exactes courtes de la forme
\begin{equation}
0\to H^i(M/\mathcal F)\to H^i(M)\to H_\kappa^{i-1}(M/\mathcal F)\to0
\end{equation}
On en déduit que $b_p(M)=\dimension H^i(M/\mathcal F)+\dimension
H_\kappa^{i-1}(M/\mathcal F)$ et donc que $m_p$=0.
Par contraposée, on obtient la première affirmation du corollaire.

Si la classe d'\'Alvarez est non nulle, alors $H^{n-1}(M/\mathcal F)=0$
(théorème \ref{coh:geod}) et donc $H_\kappa^0(M/\mathcal F)=0$
(fait \ref{coh:dual}). On peut alors déduire de la suite de Gysin
que $H^1(M/\mathcal F)\simeq H^1(M)$, et on obtient finalement que $m_1=0$.

Il reste à montrer que si $[\kappa]=0$ et $[e]\neq0$, alors on a effectivement
une petite valeur propres pour les $1$-formes. Le début de
la suite de Gysin s'écrit :
\begin{equation}
0\to H^1(M/\mathcal F)\stackrel{i}{\to}H^1(M)\stackrel{j}{\to}
H_\kappa^0(M/\mathcal F) \stackrel{\wedge e}{\to}H^2(M/\mathcal F)\to
\end{equation}
Le fait que $[\kappa]=0$ implique que $H_\kappa^0(M/\mathcal F)\simeq
H^0(M/\mathcal F)\simeq\R$, le noyau de $\wedge e$ est donc $\R$ ou $\{0\}$. 
La cohomologie $H_\kappa^0(M/\mathcal F)$
étant alors représentée par les fonctions constantes sur $M$, l'application
$\wedge e$ est nulle si et seulement si $[e]=0$. Si $[e]\neq0$, alors
le noyau de $\wedge e$ est $\{0\}$ et l'application $j$ est nulle, 
et donc $H^1(M/\mathcal F)\simeq H^1(M)$. 
Si $[\kappa]=0$ et $[e]\neq0$, on a finalement $m_1=\dimension 
H^{n-1}(M/\mathcal F)=\dimension H_\kappa^0(M/\mathcal F)=1$.
\end{demo}

\begin{remarque}\label{coh:nul}
L'absence de petite valeur propre se traduit sur la suite de Gysin par
le fait que toutes les applications $H_\kappa^{i-1}(M/\mathcal F)
\stackrel{\wedge e}{\to}H^{i+1}(M/\mathcal F)$ sont nulles. Si le flot 
est isométrique, le cas $i=1$ permet d'en déduire que la classe
d'Euler est nulle, mais le raisonnement ne se généralise pas aux 
autres flots, et on peut trouver des flots pour lesquels toutes les 
applications $\wedge e$ sont nulles
sans que la classe d'Euler soit nulle (voir paragraphe \ref{ex:euler}).
\end{remarque}

\section{Géométrie des effondrements}\label{geom}
Avant de démontrer les théorèmes \ref{intro:th2} et \ref{intro:th3},
nous allons d'abord justifier la notion d'adhérence générique:
\begin{lem}\label{geom:lem0}
Soit $M$ une variété riemannienne compacte et $\mathcal F$ un flot
riemannien sur $M$. Il existe un flot linéaire $\mathcal F_0$ 
sur un tore $T^k$ tel que toute restriction du flot $\mathcal F$
à une adhérence de feuille soit conjuguée à un flot quotient de 
$\mathcal F_0$, et que la réunion des adhérences des feuilles de 
$\mathcal F$ sur lesquelles le flot est conjugué à 
$\mathcal F_0$ forme un ouvert dense de mesure pleine de $M$.
\end{lem}
Les flots quotients de $\mathcal F_0$ sont les flots obtenus, pour tout 
$l\leq k$, en projetant linéairement $T^k$ sur un tore $T^l$ (en incluant
l'identité et les quotients d'indice fini, si $k=l$).

Cette notion de feuille générique correspond à celle de feuille régulière 
sans holonomie introduite dans \cite{am86}. Par ailleurs, on peut noter 
que les adhérences des feuilles du flot sont de dimension au plus $k$.

\begin{demo}[du lemme \protect\ref{geom:lem0}]
On va s'appuyer sur le fait que toute feuille $F$ du flot admet un 
voisinage ouvert $U$, appelé «~voisinage de Carrière~», vérifiant les 
propriétés suivantes (voir \cite{ca84a}):
\begin{itemize}
\item le voisinage $U$ est saturé par le flot;
\item il existe un difféomorphisme $U\to S^1\times T^k\times B^{n-k}$
envoyant l'adhérence de $F$ sur $S^1\times T^k\times \{0\}$;
\item le flot restreint à $U$ est conjugué au flot obtenu par la suspension
d'un difféomorphisme $T\times R$ de $T^k\times B^{n-k-1}$, où $T$
est une translation irrationnelle et $R$ une rotation de la boule $B^{n-k-1}$.
\end{itemize}
On peut décomposer la rotation $R$ en somme de rotations planes :
\begin{equation}
R=\left(\begin{array}{cccc}
R_1&\ldots&0&0\\
\vdots&\ddots&\vdots&\vdots\\
0&\ldots&R_m&0\\
0&\ldots&0&Id
\end{array}\right)
\end{equation}
où les $R_i$ sont des rotations non triviales du plan.

Supposons d'abord que $R\neq Id$.
On définit les points génériques de la boule $B^{n-k-1}$ comme étant 
les points qui ne sont fixés par aucune des rotations $R_i$. 
Si $x$ et $y$ sont deux tels points, il existe une matrice diagonale $D$ et 
une rotation $R'$ commutant toutes deux avec $R$ et telle que $y=DR'(x)$. 
Les orbites de $x$ et $y$ sous l'action de $R$ sont conjuguées
\emph{via} l'application $DR'$, donc leur orbites pour le flot le sont
aussi, ainsi que les adhérences de ces orbites. 
L'ensemble des points génériques de la boule forme un ouvert 
de mesure pleine de $B^{n-k-1}$, donc la réunion de leurs orbites pour
le flot est un ouvert de mesure pleine de $U$. On appellera feuilles génériques
(resp. adhérences génériques) de $U$ ces orbites (resp. leur adhérence). 
Si on considère un point $z$ non générique de la boule,
l'adhérence de son orbite pour le flot peut s'écrire comme le quotient
de l'adhérence générique par l'adhérence du groupe engendré par les rotation
$R'_i$ qui fixent ce point.

Si on suppose maintenant que $R=Id$, toutes les restrictions du flot
aux adhérences des feuilles de $U$ sont conjuguées au flot sur $\bar F$,
qui représente alors la feuille générique de $U$.

Par compacité, il suffit d'un nombre fini de voisinages de Carrière
pour recouvrir la variété $M$. Si deux de ces ouverts ont une 
intersection commune, ils auront en particulier une feuille en
commun qui sera une feuille générique de chacun des deux ouverts.
Par connexité, on en déduit que tous les ouverts du recouvrement
on la même adhérence générique. La réunion des feuilles génériques
formant un ouvert dense de mesure pleine sur chaque ouvert du recouvrement,
leur réunion sur l'ensemble de la variété est aussi un ouvert dense de
mesure pleine. Il en va de même pour les adhérences génériques.
\end{demo}

Quand on effondre la variété $M$ par $g_\varepsilon$, elle tend vers 
un espace métrique $X$  qui est le quotient de $M$ par les adhérences
des feuilles du flot. La distance de Gromov-Hausdorff 
$\delta(M,\mathcal F,g_\varepsilon)$ entre $(M,g_\varepsilon)$ et $X$
est alors de l'ordre du maximum du diamètre des adhérences des feuilles
pour la métrique induite. Plus précisément, on sait (voir \cite{cfg92},
en particulier le théorème 2.6) que si on note $\pi:M\to X$
la projection de $M$ sur $X$ et qu'on pose 
\begin{equation}
\delta'(\varepsilon)=\sup_{x\in X}\diam(\pi^{-1}(x), g_\varepsilon),
\end{equation}
 alors il existe une constante $c>0$ 
indépendante de $\varepsilon$ telle que
\begin{equation}\label{geom:eqd}
\frac1c\cdot\delta'(\varepsilon)\leq\delta(M,\mathcal F,g_\varepsilon)
\leq c\cdot\delta'(\varepsilon),
\end{equation}
pour tout $\varepsilon\in]0,1]$.  
En nous appuyant sur cette remarque, on va montrer que 
$\delta(M,\mathcal F,g_\varepsilon)$ est contrôlé par le diamètre
d'une adhérence générique.
\begin{lem}\label{geom:lemd}
Soit $(M,g)$ une variété riemannienne compacte munie d'un flot riemannien 
$\mathcal F$ d'adhérence générique $(T^k,\mathcal F_0)$. Si $\bar g$ est une
métrique sur $T^k$ et en notant $(g_\varepsilon)$ et $(\bar g_\varepsilon)$
les effondrements adiabatiques de $M$ et $T^k$ pour leur flot respectif,
alors il existe un constante $c>0$ ne dépendant pas de $\varepsilon$ telle
que
$$\frac1c\cdot\diam(T^k,\bar g_\varepsilon)\leq
\delta(M,\mathcal F,g_\varepsilon)\leq c\cdot\diam(T^k,\bar g_\varepsilon)$$
pour tout $\varepsilon\in]0,1]$.
\end{lem}
\begin{demo}
Pour tout $x\in X$, la métrique $g$ induit une métrique sur $\pi^{-1}(x)$. 
On peut ainsi définir une application
continue de $X$ dans l'espace des métriques du tore $T^k$: Si $x$
est l'image par $\pi$ d'une adhérence générique, l'image de $x$ est simplement
la métrique $g_x$ induite par $g$ sur cette adhérence, et
on peut prolonger cette application pas continuité sur les adhérences
non génériques, quitte à ce que la métrique $g_x$ soit dégénérée dans le cas
où l'adhérence de la feuille est de dimension strictement inférieure à $k$, 
le diamètre de $\pi^{-1}(x)$ étant dans tous les cas majoré par 
$\diam(T^k,g_x)$. 

L'espace $X$
étant compact, on peut trouver une constante $\tau>0$ telle que
pour tout $x\in X$, $g_x\leq\tau\bar g$, où $\bar g$ est
la métrique canonique de $\R^k/\Z^k$. En notant $(g_{x,\varepsilon})$
et $(\bar g_\varepsilon)$ les effondrements adiabatiques associés, on a
$g_{x,\varepsilon}\leq\tau\bar g_\varepsilon$, pour tout $\varepsilon>0$.
On en déduit que 
\begin{equation}
\diam(\pi^{-1}x,g_\varepsilon)\leq\diam(T^k,
g_{x,\varepsilon})\leq\tau\diam(T^k,\bar g_\varepsilon)
\end{equation}
pour tout $x\in X$ et tout $\varepsilon\in]0,1]$, et donc que 
\begin{equation}\label{geom:eqd1}
\delta'(\varepsilon)
\leq\tau\diam(T^k,\bar g_\varepsilon).
\end{equation}
 Si on se donne un point $x_0\in X$ qui soit la projection par $\pi$
d'une adhérence générique, on peut trouver $\tau'>0$ tel que 
$\tau'\bar g\leq g_{x_0}$. On a alors 
\begin{equation}\label{geom:eqd2}
\tau' 
\diam(T^k,\bar g_\varepsilon)\leq\delta'(\varepsilon)
\end{equation}
 pour tout $\varepsilon\in]0,1]$. On peut alors conclure
en utilisant l'encadrement~(\ref{geom:eqd}).
\end{demo}

Le lemme \ref{geom:lemd} assure que
\begin{equation}
\liminf_{\varepsilon\to0}\frac{\ln\delta(M,\mathcal F,g_\varepsilon)}
{\ln\varepsilon}=\liminf_{\varepsilon\to0}\frac{\ln\diam(T^k,
\bar g_\varepsilon)}{\ln\varepsilon}.
\end{equation}
 Pour démontrer le théorème \ref{intro:th2}, il reste donc à décrire 
la géométrie de l'effondrement d'une adhérence générique 
dans le cas où elle est de dimension 2. Le flot est entièrement déterminé 
par la donnée d'un irrationnel $\alpha\in\R$, le flot étant induit par le 
vecteur $(1,\alpha)\in\R^2$. On va
étudier le comportement asymptotique du diamètre d'un tore $T^2$ plat 
au cours de l'effondrement d'un flot linéaire:
\begin{lem}\label{geom:lem1}
Soit $T^2=\R^2/\Z^2$, $g$ le quotient sur $T^2$ de la métrique
canonique de $\R^2$, $\alpha$ un réel irrationnel et $(g_\varepsilon)$ 
l'effondrement adiabatique associé à $g$ et au flot linéaire induit
par le vecteur $(1,\alpha)$. On a alors :
\begin{equation}
\liminf_{\varepsilon\to0}\frac{\ln\diam(T^2,g_\varepsilon)}
{\ln\varepsilon}=\frac1{\mu(\alpha)}.
\end{equation}
\end{lem}
\begin{demo}
On va en fait montrer que
\begin{equation}
\limsup_{\varepsilon\to0}\frac{\ln\injrad(T^2,g_\varepsilon)}
{\ln\varepsilon}=1-\frac1{\mu(\alpha)},
\end{equation}
le lemme s'en déduit en utilisant le fait qu'il existe une constante
$C>0$ telle que
\begin{equation}\label{geom:vol}
\frac1C\cdot\Vol(T^2,g_\varepsilon)\leq\injrad(T^2,g_\varepsilon)\cdot
\diam(T^2,g_\varepsilon)\leq C\cdot\Vol(T^2,g_\varepsilon).
\end{equation}
et que $\Vol(T^2,g_\varepsilon)=\varepsilon$.

Notons $\|\cdot\|_\varepsilon$ la norme sur $\R^2$ obtenue en relevant
 la métrique $g_\varepsilon$. On peut remarquer que pour
tout $\varepsilon>0$, on a 
\begin{equation}
\injrad(T^2,g_\varepsilon)=\frac12\inf_{\underset{(p,q)\neq(0,0)}
{(p,q)\in\Z^2}} \|(p,q)\|_\varepsilon
\end{equation}
Le problème est donc d'estimer la distance pour la métrique $g_\varepsilon$
entre les points non nuls de $\Z^2$ et l'origine. On notera $D$ la droite 
vectorielle de $\R^2$ engendrée par $(1,\alpha)$, et $\theta$ l'angle entre 
la droite $D$ et l'axe des abscisses pour la métrique $g$.

On va d'abord minorer cette distance en fonction de l'exposant
d'irrationalité de $\alpha$. Soit $\nu$ un réel tel que $\nu>\mu(\alpha)$. 
L'inéquation $|\alpha-\frac pq|<\frac1{q^\nu}$ n'a alors qu'un nombre
fini de solutions, on peut donc trouver une constante $0<c(\alpha)<1$ telle que
\begin{equation}\label{geom:eqnu}
\left|\alpha-\frac pq\right|>\frac c{q^\nu}
\end{equation}
pour tout $(p,q)\in\Z^2\backslash(0,0)$. Pour tout $\varepsilon>0$ on
définit une partie $A_\varepsilon$ de $\R^2$ par:
\begin{equation}
A_\varepsilon=\left\{(x,y)\in\R^2,\ |x\alpha-y|<c\cdot\varepsilon^{1-\frac1\nu}
,\ |x|<\frac1{\varepsilon^{\frac1\nu}}\right\}.
\end{equation}
Pour tout $0<\varepsilon<1$, le seul point à coordonnées entières contenu 
dans $A_\varepsilon$
est $(0,0)$. En effet, si $(q,p)\in A_\varepsilon\cap\Z^2$, alors
$|q|<\varepsilon^{-\frac1\nu}$ et $|q\alpha-p|<
c\cdot\varepsilon^{1-\frac1\nu}$. Si $q\neq0$ ceci implique
$|q\alpha-p|<\frac c{q^{\nu-1}}$, et donc $|\alpha-\frac pq|<\frac c{q^\nu}$,
ce qui contredit (\ref{geom:eqnu}). Si $q=0$, alors 
$|p|<c\cdot\varepsilon^{1-\frac1\nu}$ ce qui implique $p=0$.

Le domaine $A_\varepsilon$ est, pour la métrique $g$, un parallélogramme
dont un coté, parallèle à la droite $D$, est de longueur $\frac2{\cos\theta}
\varepsilon^{-\frac1\nu}$, et l'autre, parallèle à l'axe
des ordonnées, est de longueur $2c\cdot\varepsilon^{1-\frac1\nu}$. Il contient
donc un rectangle $R$ dont les cotés, respectivement parallèles et orthogonaux
à $D$, sont de longueurs $\frac2{\cos\theta}
\varepsilon^{-\frac1\nu}-2c|\sin\theta|\varepsilon^{1-\frac1\nu}$ et
$2c\cdot\cos\theta\cdot\varepsilon^{1-\frac1\nu}$. Pour la métrique 
$g_\varepsilon$,
$R$ est un rectangle de cotés $\varepsilon^{1-\frac1\nu}(\frac2{\cos\theta}
-2c|\sin\theta|\varepsilon)$ et $2c\cdot\cos\theta\cdot
\varepsilon^{1-\frac1\nu}$. 
On peut donc trouver un réel $\rho>0$ ne dépendant pas de $\varepsilon$ 
tel que, pour la métrique $g_\varepsilon$, la boule
de centre $(0,0)$ et de rayon $\rho\varepsilon^{1-\frac1\nu}$
soit contenue dans $A_\varepsilon$ et donc ne contienne pas d'autre
point à coordonnées entière. On a par conséquent la minoration
\begin{equation}
\injrad(T^2,g_\varepsilon)\geq\frac\rho2\varepsilon^{1-\frac1\nu}
\end{equation}
et, en tenant compte du fait que $\ln\varepsilon$ est négatif, la majoration
\begin{equation}
\frac{\ln\injrad(T^2,g_\varepsilon)}{\ln\varepsilon}\leq1-\frac1\nu+
\frac{\ln\frac\rho2}{\ln\varepsilon}.
\end{equation}
En passant à la limite quand $\varepsilon\to0$ puis $\nu\to\mu(\alpha)$,
on obtient
\begin{equation}
\limsup_{\varepsilon\to0}\frac{\ln\injrad(T^2,g_\varepsilon)}
{\ln\varepsilon}\leq1-\frac1{\mu(\alpha)},
\end{equation}

On va maintenant majorer $\injrad(T^2,g_{\varepsilon_n})$ pour une
suite $(\varepsilon_n)$ tendant vers zéro. Soit $\nu<\mu(\alpha)$. Il existe
deux suites d'entiers $(p_n)$ et $(q_n)$ telles que $(q_n)$ soit
strictement croissante et $\left|\alpha-\frac{p_n}{q_n}\right|<\frac1{q_n^\nu}$
pour tout $n\in\N$.

Soit $P_n$ le projeté orthogonal de $(p_n,q_n)$ sur la droite $D$. 
On peut noter que $P_n$ est indépendant du choix de la métrique 
dans la famille $(g_\varepsilon)$.
Si on note $d$ et $d'$ les distances pour la métrique canonique 
de l'origine à $P_n$ et de $P_n$ à $(p_n,q_n)$, on a $\|(p_n,q_n)\|^2_\varepsilon
=\varepsilon^2d^2+d'^2$. D'un part, la distance $d'$ peut s'écrire
$d'=|q_n\alpha-p_n|\cos\theta$, donc $d'\leq q_n^{1-\nu}\cos\theta$. D'autre
part, comme $\frac{p_n}{q_n}\to \alpha$, on a la majoration
\begin{equation}
d^2\leq p_n^2+q_n^2=q_n^2\left(\frac{p_n^2}{q_n^2}+1\right)\leq q_n^2(\alpha+2).
\end{equation}
Si on pose $\varepsilon_n=q_n^{-\nu}$, on a finalement
\begin{equation}
\|(p_n,q_n)\|^2_{\varepsilon_n}\leq q_n^{2(1-\nu)}(\cos^2\theta+\alpha+2)
=\varepsilon_n^{2(1-\frac1\nu)}(\cos^2\theta+\alpha+2).
\end{equation}
On en déduit que
\begin{equation}
\limsup_{\varepsilon\to0}\frac{\ln\injrad(T^2,g_{\varepsilon})}
{\ln\varepsilon}\geq
\limsup_{n\to+\infty}\frac{\ln\injrad(T^2,g_{\varepsilon_n})}
{\ln\varepsilon_n}\geq1-\frac1\nu,
\end{equation}
et, en passant à la limite quand $\nu\to\mu(\alpha)$,
\begin{equation}
\limsup_{\varepsilon\to0}\frac{\ln\injrad(T^2,g_{\varepsilon})}
{\ln\varepsilon}\geq1-\frac1{\mu(\alpha)}.
\end{equation}
\end{demo}

\begin{remarque} Le théorème \ref{intro:th2} apporte un élément
de réponse à la question, posée par P.~Pansu dans \cite{pa84},
de savoir si dans l'exemple \ref{intro:ex} on peut mesurer l'écart
entre deux métriques obtenues en fixant $\varepsilon$ et en faisant
varier $\alpha$. On voit que si on choisit deux valeurs de $\alpha$ ayant
des exposants d'irrationalité différents, on pourra effectivement trouver
des valeurs de  $\varepsilon$ pour lesquelles la distance de Gromov-Hausdorff
entre les deux métriques est minorées.
\end{remarque}
\begin{remarque}\label{geom:rq}
 Les adhérences non génériques empêchent de contrôler
le rayon d'injectivité global en fonction de $\mu(\mathcal F)$. On
peut voir dans l'exemple \ref{intro:ex} que le rayon d'injectivité est
donné par la longueur des feuilles compactes et qu'il est donc de
l'ordre de $\varepsilon$ quel que soit $\mu(\mathcal F)$.
\end{remarque}

Une difficulté apparaît si on cherche à généraliser le lemme \ref{geom:lem1}
aux tores de dimension supérieure. On peut représenter un flot linéaire
sur $T^k$ par un vecteur $(1,\alpha)\in\R^k$ avec $\alpha\in\R^{k-1}$ et
définir un exposant d'irrationalité de $\alpha$ par 
\begin{equation}
\mu(\alpha)=\sup\{\nu,\ \|\alpha-\frac pq\|<\frac1{q^\nu} 
\textrm{ a une infinité de solutions } (p,q)\in\Z^{k-1}\times\Z\}
\end{equation}
Cet invariant permet de contrôler le rayon d'injectivité du tore lors de
l'effondrement, mais la relation (\ref{geom:vol}) se généralise par 
\begin{equation}
\frac1{C(k)}\injrad^{k-1}(T^k)\cdot\diam(T^k)\leq
\Vol(T^k)\leq C(k)\cdot\injrad(T^k)\cdot
\diam^{k-1}(T^k),
\end{equation}
ce qui ne permet pas de contrôler le diamètre en fonction du 
rayon d'injectivité aussi précisément qu'en dimension 2. Cependant,
dans le cas où $\alpha$ est difficilement approchable, on a montré dans
\cite{ja04} le résultat suivant :
\begin{proposition}
Soit $\alpha\in\R^{k-1}$ un vecteur difficilement approchable, $\mathcal F$ le
flot linéaire induit sur $T^k$ par le vecteur $(1,\alpha)$, $g$
une métrique sur $T^k$ et $(g_\varepsilon)$ l'effondrement 
adiabatique associé à $\mathcal F$ et $g$. Il existe une constante 
$c(g,\alpha)>0$ telle que
\begin{equation}
\frac1c\cdot\varepsilon^{\frac1k}\leq\injrad(T^k,g_\varepsilon)\leq
\diam(T^k,g_\varepsilon)\leq c\cdot\varepsilon^{\frac1k}.
\end{equation}
\end{proposition}

 Le lemme \ref{geom:lemd} permet d'en déduire directement le théorème 
\ref{intro:th3}.

\section{Exemples}\label{ex}
\subsection{Flots isométriques}
\subsubsection{Effondrements de sphères}\label{ex:spheres}
Les sphères ayant un gros groupe d'isométrie, elles permettent de construire
facilement des exemples de flots isométriques. En outre, on 
obtient aisément des flots pour lesquels l'espace limite de l'effondrement
est singulier. On peut généraliser l'exemple \ref{intro:ex} de la
manière suivante :

\begin{exemple}
On considère la sphère 
\begin{equation}
S^{2n-1}=\{(a_1,\ldots,a_n)\in\C^n,\ |a_1|^2+\ldots+|a_n|^2=1\}
\end{equation}
sur laquelle le tore $T^n$ agit isométriquement par
\begin{equation}
(\theta_1,\ldots,\theta_n)\cdot(a_1,\ldots,a_n)=
(e^{i\theta_1}a_1,\ldots,e^{i\theta_n}a_n).
\end{equation}
Un sous-groupe à un paramètre du tore $T^n$ induit un flot isométrique
non singulier sur $S^{2n-1}$. Si ce sous-groupe est dense dans $T^n$, 
l'espace limite
de l'effondrement est $X=S^{2n-1}/T^n$, qu'on peut identifier à un domaine
de la sphère $S^{n-1}$:
\begin{equation}
X=\{(x_1,\ldots,x_n)\in\R^n,\ |x_1|^2+\ldots+|x_n|^2=1\textrm{ et }
x_1,\ldots,x_n\geq0\}.
\end{equation}
Cet espace est difféomorphe à un simplexe de dimension $n-1$.

Le même argument que dans l'exemple \ref{intro:ex} assure que ces
effondrements produisent des petites valeurs propres. La classe d'Euler
de ces flots est donc non nulle, selon le corollaire \ref{intro:cor1}.
\end{exemple}

\begin{remarque} Comme les sphères de dimension paire ont une caractéristique 
d'Euler non nulle, non seulement elles n'admettent pas de flot non singulier
mais elles n'admettent pas non plus d'effondrement à courbure bornée.
\end{remarque}

\begin{remarque}
La suite de Gysin permet de calculer la cohomologie basique du flot,
on a $H^p(M/\mathcal F)=\R$ quand $p$ est pair et 
$H^p(M/\mathcal F)=\{0\}$ quand $p$ est impair. On a donc des petites
valeurs propres en tout degré autre que 0 et $n$.
\end{remarque}

\subsubsection{Prescription de l'exposant d'irrationalité}\label{ex:ex}
Dans l'exemple \ref{intro:ex}, comme sur n'importe quelle variété
dont le groupe d'isométrie contient un tore $T^2$, 
on peut construire
pour tout $\alpha\in\R\backslash\Q$  un flot dont le vecteur directeur
sur une adhérence générique est $(1,\alpha)$. Pour tout $\mu\in[2,+\infty]$,
on peut construire un $\alpha$ dont l'exposant d'irrationalité est
$\mu$ en utilisant le fait que l'exposant d'irrationalité peut se calculer
à l'aide du développement en fraction continue de $\alpha$:
\begin{theo}[\cite{so04}]
Soit $\alpha$ un irrationnel, $[a_0,a_1,a_2,\ldots]$ son développement en 
fraction continue et $\displaystyle\frac{p_n}{q_n}$ les réduites
successives de ce développement. On a alors
$$\mu(\alpha)=1+\limsup_{n\to+\infty}\frac{\ln q_{n+1}}{\ln q_n}
=2+\limsup_{n\to+\infty}\frac{\ln a_{n+1}}{\ln q_n}.$$
\end{theo}
Sachant qu'on a aussi $q_{n+1}=a_{n+1}q_n+q_{n-1}$, on voit qu'on peut 
aisément construire une suite $(a_n)$ ayant l'exposant d'irrationalité
souhaité, par exemple en posant la relation de récurrence $a_{n+1}=[q_n^\mu]$
si $\mu$ est fini. Le cas $\mu=+\infty$ correspond aux nombres
de Liouville.

 Si le groupe d'isométrie contient un tore $T^k$, on peut construire
un flot linéaire sur ce tore ---~et donc un flot isométrique sur la
variété~--- dont la direction est difficilement approchable, en utilisant 
le fait que si $(1,\alpha_1,\ldots,\alpha_m)$ est
une base d'un corps de nombres réel, alors $(\alpha_1,\ldots,\alpha_m)$
est difficilement approchable (\cite{sc80}, théorème 4A, p.42).

\subsubsection{Fibrés principaux en tores}
Si on considère un fibré $T^2\hookrightarrow M\to N$ principal, un
vecteur $(1,\alpha)$ sur $T^2$ induit un flot isométrique dont
toutes les feuilles sont génériques. En suivant la construction donnée
dans \cite{ja04}, on peut obtenir un effondrement pour lequel le
rayon d'injectivité global est du même ordre que celui d'une fibre 
quelconque. On peut donc construire, pour tout $\mu\geq2$, un flot pour
l'effondrement duquel on a
\begin{equation}
\liminf_{\varepsilon\to0}\frac{\ln\lambda_{p,i}(M,g_\varepsilon)}
{\ln\injrad(M,g_\varepsilon)}=\frac{2\mu}{\mu-1},
\end{equation}
Pour tout $i\leq m_p$.

\subsection{Flots non isométriques}
\subsubsection{Fibrés en tores sur le cercle}
On peut construire un grand nombre d'exemples de flots non isométriques
en prenant comme support des fibrés en tores sur le cercle. Si on
se donne un matrice $A\in\SL_k(\Z)$, cette matrice induit un difféomorphisme
linéaire du tore $T^k$ et on peut construire un fibré 
$T^k\hookrightarrow M\to S^1$ par suspension de ce difféomorphisme.
Chaque vecteur propre de $A$ induit sur $T^k$ un flot linéaire qui s'étend
en un flot global sur $M$. Si on munit $M$ d'une structure de solvariété
et d'une métrique invariante comme dans \cite{ja03}, on peut aisément
représenter le flot considéré par un champ de vecteur invariant et
vérifier que la métrique est quasi-fibrée pour ce flot (voir les
exemples donnés dans \cite{rp01b}). 

Pour illustrer le corollaire \ref{intro:cor1}, on va donner deux exemples
illustrant chacun des cas $[e]=0$ et $[e]\neq0$.

\begin{exemple}\label{ex:ex1}
On pose $A=\smat{2&1\\1&1}$. Le flot correspondant est l'exemple le
plus simple de flot non isométrique ---~il représente l'unique classe
de flots non isométriques dans la classification des flots riemanniens
établie en dimension~3 par Y.~Carrière dans \cite{ca84a}~---, et sa classe
d'Euler est nulle (cf. \cite{rp01b}). L'effondrement du flot ne produit
donc aucune petite valeur propre.
\end{exemple}

\begin{exemple}\label{ex:ex2}
 On pose
\begin{equation}
A=\left(\begin{array}{cccc}
2&1&0&1\\ 1&1&0&0\\ 0&0&2&1\\ 0&0&1&1
\end{array}\right).
\end{equation}
D'après la théorie de J.~Lott (\cite{lo02}), les effondrements du
fibré obtenu ne produisent jamais de petites valeurs propres pour les
$1$-formes différentielles, mais on a montré dans \cite{ja03} que
l'effondrement du flot produisait des petites valeurs propres sur
les $2$-formes. On peut en déduire que $[\kappa]\neq0$ et $[e]\neq0$.
\end{exemple}
 Dans \cite{rp01b} est présenté un exemple de variété de dimension 6
admettant un flot riemannien vérifiant $[\kappa]\neq0$ et $[e]\neq0$. Dans
l'exemple \ref{ex:ex2} ci-dessus, la variété est de dimension 5. On ne 
peut pas trouver de tels exemples en dimension 3 selon la classification 
de Y.~Carrière. On va voir au paragraphe suivant qu'il en existe en
dimension 4, mais que leur effondrement ne peut pas produire de
petites valeurs propres.

\subsubsection{Existence de petites valeurs propres et 
annulation de la classe d'Euler}\label{ex:euler}
Nous allons construire dans ce paragraphe  un exemple de flot dont
la classe d'Euler est non nulle mais dont l'effondrement adiabatique ne
produit pas de petite valeur propre. Pour cela, nous commencerons par
établir deux lemmes. Le premier fournit un critère concernant l'annulation
des classes d'Euler et d'Álvarez :
\begin{lem}\label{ex:lem1}
Soit $\mathcal F$ un flot riemannien sur une variété $M$, et $N$ une
sous-variété de $M$ stable par $\mathcal F$. Si la classe d'Álvarez
(resp. d'Euler) de $\mathcal F$ est nulle, alors la classe d'Álvarez
(resp. d'Euler) du flot restreint à $N$ est nulle aussi.
\end{lem}
\begin{demo}
Le cas de la classe d'Álvarez est évident quand on se souvient que
$[\kappa]=0$ signifie que le flot est isométrique: si une métrique sur 
$M$ est invariante par le flot, sa restriction à $N$ sera aussi invariante
sous l'action de la restriction du flot.

Si $\chi$ est la forme caractéristique du flot sur $M$, la forme
d'Euler est définie par $\de\chi=\kappa\wedge\chi+e$, et $[e]=0$ signifie
que $e$ est une forme exacte pour la différentielle tordue $\de_{-\kappa}$,
c'est-à-dire qu'il existe une forme basique $\beta$ telle que $e=\de\beta
+\kappa\wedge\beta$. On a donc $\de\chi=\kappa\wedge\chi+\de\beta
+\kappa\wedge\beta$. Si on note $i$ l'injection $i:N\to M$, on peut
en déduire $\de(i^*\chi)=(i^*\kappa)\wedge(i^*\chi)+\de(i^*\beta)
+(i^*\kappa)\wedge(i^*\beta)$. On peut vérifier que $i^*\chi$ est bien
la forme caractéristique du flot restreint à $N$, et que $i^*\kappa$ 
et $i^*\beta$ sont des formes basiques pour ce flot. On en déduit que
$i^*\kappa$ est la forme de courbure et $\de(i^*\beta)+
(i^*\kappa)\wedge(i^*\beta)=\de_{-i^*\kappa}i^*\beta$ la forme d'Euler
du flot sur $N$, et que la classe d'Euler correspondante est nulle.
\end{demo}

Le second lemme concerne le spectre du laplacien en dimension 4 :
\begin{lem}
L'effondrement adiabatique d'un flot riemannien non isométrique 
sur une variété de dimension~4 ne produit pas de petites valeurs propres.
\end{lem}
\begin{demo}
La théorie de Hodge assure qu'en dimension 4, le spectre du laplacien
de Hodge-de~Rham est entièrement déterminé par sa restriction aux $1$-formes
différentielles. Or, on a vu que l'effondrement d'un flot non isométrique 
ne produit pas de petite valeur propre sur les $1$-formes.

On peut aussi remarquer que pour un flot non isométrique en dimension~4, 
on a $H_\kappa^0(M/\mathcal F)=H^3(M/\mathcal F)=\{0\}$, et donc
toutes les applications $H_\kappa^{i-1}
(M/\mathcal F)\stackrel{\wedge e}{\to}H^{i+1}(M/\mathcal F)$ dans la suite
de Gysin du flot sont nulles, ce qui garantit qu'il n'y a pas de petites
valeurs propres (cf. remarque \ref{coh:nul}).  
\end{demo}

En vertu du lemme précédent, il nous suffit de construire un flot 
riemannien non isométrique sur une variété de dimension~4 dont la classe 
d'Euler est non nulle. L'exemple que nous allons considérer appartient
à une famille de flots non isométriques exhibée dans \cite{am86}, nous
allons en donner une construction qui permet de facilement mettre en
évidence le fait que $[\kappa]\neq0$ et $[e]\neq0$ en utilisant le lemme
\ref{ex:lem1}.

On commence par considérer le flot de l'exemple \ref{intro:ex} sur la sphère
(le coefficient $\alpha$ déterminant la direction du flot sera fixé
ultérieurement).
En prenant le produit de la sphère $S^3$ par un cercle, on obtient un 
flot sur $S^3\times S^1$ dont l'adhérence générique est de dimension~2 et 
dont l'espace limite de l'effondrement est $[0,\frac\pi2]\times S^1$. 
On peut noter qu'au dessus de $]0,\frac\pi2[\times S^1$, le flot se
décompose en un produit trivial entre l'adhérence générique et 
$]0,\frac\pi2[\times S^1$.

\begin{figure}[h]
\begin{center}
\begin{picture}(0,0)%
\includegraphics{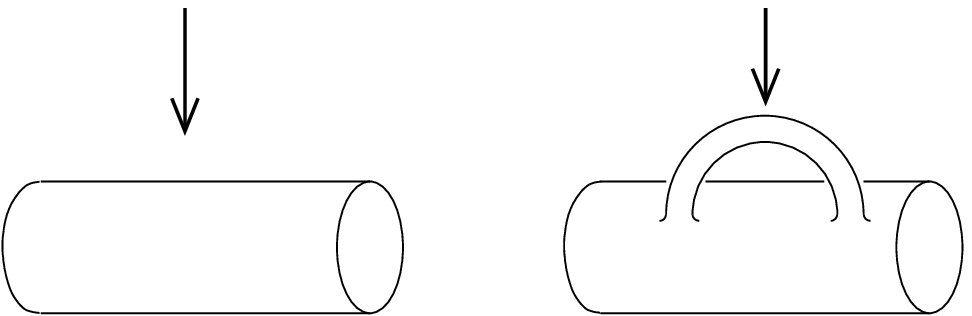}%
\end{picture}%
\setlength{\unitlength}{4144sp}%
\begingroup\makeatletter\ifx\SetFigFontNFSS\undefined%
\gdef\SetFigFontNFSS#1#2#3#4#5{%
  \reset@font\fontsize{#1}{#2pt}%
  \fontfamily{#3}\fontseries{#4}\fontshape{#5}%
  \selectfont}%
\fi\endgroup%
\begin{picture}(4409,1689)(1315,-2548)
\put(1846,-1006){\makebox(0,0)[lb]{\smash{{\SetFigFontNFSS{12}{14.4}{\rmdefault}{\mddefault}{\updefault}{\color[rgb]{0,0,0}$S^3\times S^1$}%
}}}}
\put(2251,-1411){\makebox(0,0)[lb]{\smash{{\SetFigFontNFSS{12}{14.4}{\rmdefault}{\mddefault}{\updefault}{\color[rgb]{0,0,0}$\pi$}%
}}}}
\put(4906,-1366){\makebox(0,0)[lb]{\smash{{\SetFigFontNFSS{12}{14.4}{\rmdefault}{\mddefault}{\updefault}{\color[rgb]{0,0,0}$\pi$}%
}}}}
\put(4726,-1006){\makebox(0,0)[lb]{\smash{{\SetFigFontNFSS{12}{14.4}{\rmdefault}{\mddefault}{\updefault}{\color[rgb]{0,0,0}$M$}%
}}}}
\end{picture}%
\end{center}
\caption{\label{ex:fig1}}
\end{figure}

On modifie maintenant la variété de manière à ajouter une anse à l'espace
limite: le nouvel espace limite $X$ est obtenu en enlevant deux disques $D_0$
et $D_1$ à
$[0,\frac\pi2]\times S^1$ et en recollant le long des bords de ces deux
disques un cylindre $S^1\times[0,1]$. La variété $M$ sur laquelle sera
défini le nouveau flot est obtenu en enlevant les deux parties de la forme
$D_i\times T^2$ au dessus des deux disques de l'espace limite, et en
recollant le long du bord le produit $S^1\times[0,1]\times T^2$. On doit
préciser la manière d'effectuer le recollement car la topologie
de la variété $M$ en dépend: on se donne une matrice $A\in\SL_2(\Z)$ et
on identifie les bords de $D_i\times T^2$ et $S^1\times[0,1]\times T^2$
par
\begin{equation}
\left\{\begin{array}{ccc}
\partial D_0\times T^2 & \to & S^1\times\{0\}\times T^2 \\
(x,y) & \to & (x,y) \end{array}\right.
\end{equation}
et
\begin{equation}
\left\{\begin{array}{ccc}
\partial D_1\times T^2 & \to & S^1\times\{1\}\times T^2 \\
(x,y) & \to & (x,Ay) \end{array}\right. ,
\end{equation}
l'identification entre $\partial D_i$ et $S^1$ provenant de l'adjonction
d'anse sur l'espace limite.

 On construit le flot riemannien sur $M$ en munissant
$S^1\times[0,1]\times T^2$ du flot obtenu par produit de l'adhérence 
générique du flot sur $S^3\times S^1$ par $S^1\times[0,1]$. Pour assurer
la compatibilité du flot avec le recollement ci-dessus, on choisit la
matrice $A$ telle qu'elle admette une valeur propre irrationnelle
---~par exemple la matrice $A$ de l'exemple \ref{ex:ex1}~---, et
on prend comme direction du flot sur $S^3$ la direction propre correspondante.
On obtient bien ainsi un flot linéaire dense dans $T^2$ qui est invariant par 
$A$.

\begin{figure}[h]
\begin{center}
\begin{picture}(0,0)%
\includegraphics{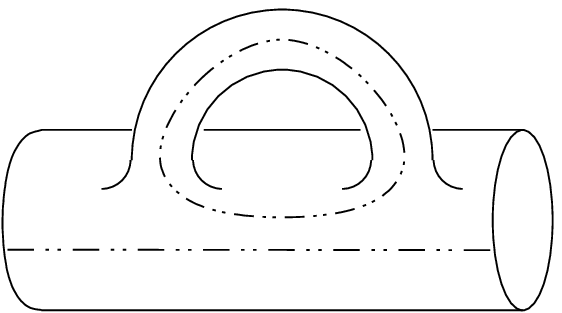}%
\end{picture}%
\setlength{\unitlength}{4144sp}%
\begingroup\makeatletter\ifx\SetFigFontNFSS\undefined%
\gdef\SetFigFontNFSS#1#2#3#4#5{%
  \reset@font\fontsize{#1}{#2pt}%
  \fontfamily{#3}\fontseries{#4}\fontshape{#5}%
  \selectfont}%
\fi\endgroup%
\begin{picture}(2529,1397)(3104,-1880)
\put(4856,-1758){\makebox(0,0)[lb]{\smash{{\SetFigFontNFSS{10}{12.0}{\rmdefault}{\mddefault}{\updefault}{\color[rgb]{0,0,0}$X_1$}%
}}}}
\put(4911,-1401){\makebox(0,0)[lb]{\smash{{\SetFigFontNFSS{10}{12.0}{\rmdefault}{\mddefault}{\updefault}{\color[rgb]{0,0,0}$X_2$}%
}}}}
\end{picture}%
\end{center}
\caption{\label{ex:fig2}}
\end{figure}
Il reste à montrer que les classes d'Álvarez et d'Euler du flot sont
non nulles. On va pour cela définir deux sous-variétés $N_1$ et $N_2$ de $M$, 
stable par le flot, et appliquer le lemme \ref{ex:lem1}.
Pour $N_1$, on choisit une sphère $S^3\times\{x\}$ dans $S^3\times S^1$
qui ne rencontre pas les $D_i\times T^2$, et qui est donc préservée dans $M$.
L'espace limite de l'effondrement de $N_1$ est un segment $X_1$ dans $X$.
 On construit $N_2$ en se donnant un cercle $X_2$ dans $X$ qui passe par 
l'anse (voir figure \ref{ex:fig2}) et en définissant $N_2$ comme étant le
fibré en tore $T^2$ au dessus de $X_2$.

La sous-variété $N_1$ est une sphère donc le flot restreint à $N_1$
a une classe d'Euler non nulle, et $N_2$ est la suspension du difféomorphisme
linéaire de $T^2$ induit par la matrice $A$, donc le flot restreint à $N_2$ 
est non isométrique d'après la classification de \cite{ca84a}. En vertu du 
lemme \ref{ex:lem1}, les classes d'Álvarez et d'Euler du flot sur $M$ sont
toutes les deux non nulles.

\begin{remarque}
Le procédé d'adjonction d'anse sur l'espace limite peut aisément se 
généraliser à d'autres flots, par exemple ceux du paragraphe \ref{ex:spheres},
pour construire des exemples de flots non isométriques.
\end{remarque}

\subsubsection{Construction de flots difficilement 
approchables}\label{ex:diff}
 On va présenter ici une construction de flots riemanniens 
non isométriques dont la direction est difficilement approchable.
Notons qu'une construction semblable apparaît aussi dans \cite{gh83}
pour d'autres de ses propriétés. 

 On considère une variété compacte $B$ telle que $b_1(B)\neq0$, 
et on note $\tilde B$ son revêtement universel.
L'idée est de se donner un flot linéaire sur $T^k$ dont la direction
est difficilement approchable, et de construire un 
fibré en tore $T^k$ sur $B$ dont le groupe de structure est formé de 
difféomorphismes linéaires de $T^k$ préservant le flot orienté sur le tore.

On se donne une extension algébrique réelle
$\mathbb K$ de $\Q$ de degré $k$, on note $r$ le nombre de plongements 
réels de $\mathbb K$ et $s$ son nombre de paires de plongements complexes
conjugués. On montre en appendice que si 
$v\in\R^k$ est un vecteur dont les composantes forment une base d'entiers 
de $\mathbb K$ ---~c'est-à-dire une base de $\mathbb K$ sur $\Q$ qui est 
aussi une base de l'anneau des entiers de $\mathbb K$,
vu comme $\Z$-module~---, alors $v$ est difficilement approchable
et le groupe $G_v$ défini par
\begin{equation}
G_v=\{A\in\SL_k(\Z),\ Av=\lambda v\text{ avec }\lambda>0\}.
\end{equation}
est un groupe abélien libre de rang $r+s-1$, et ce rang peut prendre toutes
les valeurs entre $[\frac{k+1}2]$ et $k-1$ selon le choix du corps 
$\mathbb K$ (voir théorème~\ref{ap:th3} et remarque~\ref{ap:rq3}).

On obtient un fibré $M$ en tore sur $B$ en se donnant un morphisme de
groupe $\varphi:\pi_1(B)\to G_v$ et en définissant $M$ comme
le quotient de $\tilde B\times T^k$ pour l'action de $\pi_1(B)$
définie par $\gamma\cdot(x,y)=(\gamma\cdot x,\varphi(\gamma)\cdot y)$. 
Comme la direction
de $v$ est invariante par $G_v$, le flot induit par $v$ sur 
$\tilde B\times T^k$ passe au quotient sur $M$ en un flot difficilement
approchable.

 Si l'image de $\varphi$ est non triviale, le flot n'est pas isométrique:
en effet, si le flot était isométrique, son action sur $M$ s'étendrait
en une action de son adhérence dans le groupe d'isométrie de $M$, 
cette adhérence  étant compacte et isomorphe au tore $T^k$.
le fibré serait alors principal et son groupe de structure pourrait
se réduire à des translation sur le tore $T^k$.

\begin{remarque}
L'image de $\varphi$ dans $G_v$ est un quotient abélien de $\pi_1(B)$,
donc son rang est majoré par $b_1(B)$. L'hypothèse $b_1(B)\neq0$ assure
que l'on peut choisir l'application $\varphi$ telle que son image soit non 
triviale.
\end{remarque}

\begin{remarque}
On peut appliquer à la variété obtenue un procédé chirurgical semblable
à celui pratiqué au paragraphe \ref{ex:euler} pour construire un
flot dont la classe d'Euler est non nulle. 
\end{remarque}

\begin{remarque} Les éléments du groupe $G_v$ ont $r$ directions propres
simples en commun. On peut donc construire $r$ flots difficilement 
approchables distincts sur la variété $M$ obtenue.
\end{remarque}

\newpage

\renewcommand{\thesection}{Appendice}
\section{Sur les propriétés arithmétiques des flots riemanniens}
\renewcommand{\theequation}{A.\arabic{equation}}

À un flot riemannien dont l'adhérence générique est de dimension $k$,
on peut associer un sous-groupe de $\SL_k(\Z)$ de la manière suivante:
si $v$ est un vecteur de $\R^k$, on définit le groupe $G_v$ par
\begin{equation}
G_v=\{A\in\SL_k(\Z),\ Av=\lambda v\text{ avec }\lambda>0\}.
\end{equation}
Si le flot linéaire sur l'adhérence générique d'une flot riemannien 
$\mathcal F$ est induit par un vecteur $v$, on définit $G_\mathcal F$ 
comme étant égal au groupe $G_v$, c'est le groupe des difféomorphismes
linéaires du tore qui préservent le flot linéaire orienté induit par $v$
sur $T^k$. Le vecteur $v$ dépendant d'un choix
de base dans $\Z^k$, le groupe $G_\mathcal F$ est défini à conjugaison 
près dans $\SL_k(\Z)$. 

Les groupes $G_v$ interviennent dans la classification des feuilletages
totalement géodésiques de dimension 1 établie dans \cite{gh83}, où
É.~Ghys montre en particulier que ces groupes sont
nécessairement abéliens libres et de rang inférieur ou égal à $k-1$. 
Le groupe $G_\mathcal F$ possède aussi la propriété
suivante, qui a été remarquée dans des situations particulières
(\cite{gh83}, \cite{am86}) mais dont nous n'avons pas trouvé de démonstration 
générale dans la littérature :
\begin{theo}\label{ap:th1} 
Si le groupe $G_\mathcal F$ est trivial, alors le flot $\mathcal F$ est 
isométrique.
\end{theo}
\begin{demo}
On commence par appliquer le procédé de désingularisation du flot
décrit par P.~Molino dans \cite{mo77} (voir aussi \cite{ca84a} pour
son application aux flots): si $\mathcal F$ est un flot
riemannien sur une variété $M^n$ et qu'on note $\tilde M$ la variété
de dimension $n(n-1)$ définie comme étant le fibré des repères orthonormés
directs transverses au flot, le flot $\mathcal F$ induit un flot riemannien
$\tilde{\mathcal F}$ sur $\tilde M$, dont toutes les feuilles sont génériques
et qui a la propriété d'avoir la même 
adhérence générique que $\mathcal F$. De plus, la variété $\tilde M$ possède 
une structure de fibré sur une variété $B$ dont les fibres sont les adhérences 
$T^k$ des feuilles du flot.

 Le groupe de structure du fibré $\tilde M\to B$ est composé des
difféomorphismes du tore qui préservent le flot linéaire sur l'adhérence
générique. D'après le lemme 4.7 de \cite{ep84}, on peut réduire ce
groupe de structure aux difféomorphismes affines du tore qui préservent
le flot, leur partie linéaire étant alors nécessairement contenue dans
$G_{\tilde{\mathcal F}}=G_\mathcal F$. Si le groupe $G_\mathcal F$ 
est trivial, la fibration $\tilde M\to B$ est principale. Le flot
$\tilde{\mathcal F}$ est donc isométrique, car on peut munir $\tilde M$
d'une métrique invariante pour l'action de $T^k$ qui sera nécessairement
invariante par le flot.

On conclut en utilisant le fait (\cite{ms85}, proposition I.1) que 
$\mathcal F$ est isométrique si et seulement si $\tilde{\mathcal F}$ 
est isométrique. 
\end{demo}

On va étudier plus précisément les groupes $G_v$ afin de pouvoir 
déterminer si un flot est difficilement approchable et pour en 
construire des exemples. On s'intéresse au cas où le flot linéaire
induit par $v$ sur le tore $T^k$ est à orbites denses, c'est-à-dire
que les composantes de $v$ sont linéairement indépendantes sur $\Q$. 
Si on note $v=(v_1,\ldots,v_k)$, on peut donc associer à $v$ un espace
vectoriel sur $\Q$ de dimension $k$, noté $E_v$, par $E_v=\mathrm{vect}_\Q
(1,\frac{v_2}{v_1},\ldots,\frac{v_k}{v_1})$.

\begin{theo}\label{ap:th2}
 Si $A$ est un élément de $G_v$, alors la matrice $A$ est 
diagonalisable sur $\C$, ses valeurs propres ont toutes la même 
multiplicité.

Si de plus le vecteur $v$ est vecteur simple de $A$, alors $E_v$ est
un corps de nombres et $v$ est difficilement approchable.
\end{theo}
\begin{demo}
On commence par remarquer, comme dans la démonstration du lemme 5.3 de
\cite{gh83}, que l'indépendance
linéaire des composantes de $v$ implique que si $M$ est une matrice
à coefficients entiers telle que $Mv=0$ alors $M=0$. C'est en particulier
vrai pour les matrices de la forme $P(A)$, où $P$ est un polynôme à 
coefficients entiers et $A$ un élément de $G_v$. Soit
$\lambda$ le réel tel que $Av=\lambda v$. Si $P$ et un polynôme annulateur
de $\lambda$, alors $P(A)v=P(\lambda)v=0$, donc $P$ est un polynôme 
annulateur de $A$. Réciproquement, si $P(A)=0$, alors 
$P(\lambda)=0$. Le polynôme minimal de $\lambda$ est donc
aussi polynôme minimal de $A$. Comme les racines
du polynôme minimal d'un réel algébrique sont nécessairement simples,
la matrice $A$ est diagonalisable sur $\C$. De plus, comme le
polynôme caractéristique de $A$ est à coefficient entier, il reste
invariant par l'action du groupe de Galois du polynôme minimal de
$\lambda$, donc toutes ses racines ont la même multiplicité.

La suite de la démonstration s'appuie sur une autre remarque faite par 
É.~Ghys dans \cite{gh83}, à savoir que l'application $\varphi:G_v
\to\R^*$ qui à une matrice $A\in G_v$ associe sa valeur propre
associée à $v$ est un morphisme de groupe dont l'image est
contenue dans $E_v$. En reprenant les notations précédentes,
on a donc $\lambda\in E_v$, mais aussi $\lambda^i\in E_v$
pour tout $i\in\Z$, et par conséquent $\Q(\lambda)\subset E_v$.
Si $v$ est vecteur simple d'une matrice $A\in G_v$, alors
son polynôme minimal est de degré $k$, donc $[\Q(\lambda):\Q]=
k=\dimension E_v$, c'est-à-dire que $\Q(\lambda)=E_v$.
L'espace vectoriel $E_v$ est alors un corps de nombres réel, dont
une base est formé par la famille
$(1,\frac{v_2}{v_1},\ldots,\frac{v_k}{v_1})$. Par conséquent le vecteur 
$v$ est difficilement approchable (\cite{sc80}, théorème 4A, p.42).
\end{demo}

\begin{remarque}\label{ap:rq1}
Le théorème \ref{ap:th2} assure que les multiplicités des valeurs propres des 
éléments de $G_v$ sont nécessairement des diviseurs de $k$. 
Si on considère un flot riemannien $\mathcal F$ dont la dimension $k$ de 
l'adhérence générique est un nombre premier
et que le groupe $G_\mathcal F$ n'est pas trivial ---~par exemple si le
flot $\mathcal F$ n'est pas isométrique, selon le théorème \ref{ap:th1}~---
il contient un élément dont
les valeurs propres sont simples, et donc le flot est difficilement
approchable.
\end{remarque}

\begin{remarque}
Dans \cite{gh83}, É.~Ghys introduit, pour la donnée d'une forme 
$\omega\in\R^{n*}$ telle que $\omega(v)\neq0$, le groupe $G_{v,\omega}$
des éléments $A\in G_v$ tels que $\omega$ soit vecteur propre de $\trans A$.
On peut déduire du théorème \ref{ap:th2} que si $A\in G_v$, il existe 
nécessairement une forme propre $\omega$ de $\trans A$ telle que 
$\omega(v)\neq0$, et si $v$ est vecteur simple de $A$ cette forme 
est unique, à une constante multiplicative près. Comme les éléments de
$G_v$ commutent, cette forme sera aussi forme propre des transposées des
autres éléments de $G_v$. Sous l'hypothèse que $v$ est vecteur simple
d'un élément de $A$, on a donc $G_{v,\omega}=G_v$.
\end{remarque}

 On va donner maintenant une description plus globale des groupes $G_v$:
on peut en effet les classifier en les interprétant comme une 
représentation d'un
sous-groupe du groupe des unités d'un corps de nombres réel.
On obtient en particulier une méthode pour construire des groupes
$G_v$ non triviaux, dans le but de construire des flots riemanniens
non isométriques difficilement approchables (cf. paragraphe \ref{ex:diff}).

Précisons quelques notations: si on considère un corps de nombres réel 
$\mathbb K$ de rang $k$, on note $\mathcal O_\mathbb K$
son anneau des entiers, $\mathcal U_\mathbb K=\mathcal O_\mathbb K^*$ 
son groupe des unités et $\mathcal U_\mathbb K^+$ le sous-groupe d'indice 
fini de $\mathcal U_\mathbb K$ formé des éléments positifs dont la norme
(par rapport à l'extension $\mathbb K/\Q$) est égale à 1. On peut alors
écrire:
\begin{theo}\label{ap:th3}
Soit $\mathbb K$ un corps de nombres réel et $\mathcal B=(b_1,\ldots,b_k)$
une base d'entiers de $\mathbb K$. Il existe un morphisme 
injectif $\mathcal U_\mathbb K^+\hookrightarrow\SL_k(\Z)$ dont l'image,
que noterons $G_{\mathbb K,\mathcal B}$, est 
le groupe $G_v$ pour le vecteur $v=(b_1,\ldots,b_k)\in\R^k$, qui est
difficilement approchable.

Réciproquement, Si $v\in\R^k$ est un vecteur dont les composantes sont
linéairement indépendantes sur $\Q$, il existe un entier $d$ qui divise $k$,
une extension algébrique réelle $\mathbb K$ de $\Q$ de degré $d$, une
base d'entiers $\mathcal B$ de $\mathbb K$ et 
une matrice $P\in\GL_k(\Q)$ tels que $G_v=P^{-1}G'P\cap\SL_k(\Z)$, 
où $G'$ est le groupe
$$G'=\left\{\left(\begin{array}{ccc}A&0&0\\0&\ddots&0\\0&0&A\end{array}
\right),\ A\in G_{\mathbb K,\mathcal B}\right\}.$$
\end{theo}
\begin{remarque}
Le rang du groupe $G_v$ est majoré par le rang du groupe 
$G_{\mathbb K,\mathcal B}$ qui est égal au rang de $\mathcal U_\mathbb K$.
Si $r$ et $s$ sont respectivement le nombre de plongement réels 
de $\mathbb K$ et le nombre de paires de plongements complexes conjugués,
on a $k=r+2s$ et le théorème des unités de Dirichlet affirme que
le rang de $\mathcal U_\mathbb K$ est majoré par $r+s-1$. On retrouve
donc le fait que le rang de $G_v$ est au plus $k-1$. Remarquons aussi que
si les valeurs propres de $A\in G_v$ sont simples, alors $r$ est le nombre de
valeurs propres réelles de $A$ et $2s$ son nombre de valeurs propres complexes.
La seule connaissance de $A$ permet donc de majorer le rang de $G_v$.
\end{remarque}

\begin{remarque}\label{ap:rq3}
Comme $k=r+2s$ et $r\geq1$ (il existe au moins un plongement
réel de $\mathbb K$, à savoir l'identité) le rang du groupe 
$G_{\mathbb K,\mathcal B}$ est compris entre $[\frac{k+1}2]$ et $k-1$.
\end{remarque}

\begin{remarque}
La matrice $P$ du théorème \ref{ap:th3} et son inverse ne sont pas 
nécessairement entières,
donc les éléments de $P^{-1}G'P$ ne sont \emph{a priori} pas contenu
dans $\SL_k(\Z)$. Cependant, une partie du groupe $G'$ peut
rester entière à travers cette conjugaison. Par exemple, considérons
les matrices
$$A=\left(\begin{array}{ccc}0&0&1\\1&0&6\\0&1&0\end{array}\right),\ 
B=\left(\begin{array}{ccc}-3&-2&-6\\-6&-15&-38\\-2&-6&-15\end{array}\right)
\textrm{ et }P=\left(\begin{array}{ccc}2&0&0\\0&1&0\\0&0&1\end{array}\right).$$
$A$ et $B$ sont des matrices diagonalisables qui commutent, donc
sont éléments d'un $G_v$, où $v$ est un vecteur propre commun à $A$ et $B$.
La matrice $P^{-1}AP$ n'est pas entière, mais $P^{-1}BP$ l'est.
\end{remarque}

\begin{remarque}
Déterminer une famille de générateurs du groupe des unités d'un corps
de nombre, et donc des groupes $G_v$, est un problème difficile. On
connait un algorithme pour les corps quadratiques,  mais il n'existe rien 
de tel dès que le degré du corps dépasse 2.
\end{remarque}

\begin{demo}[du théorème \protect\ref{ap:th3}]
Soit $\mathbb K$ un corps de nombres réel, $\mathcal B=(b_1,\ldots,b_k)$
une base d'entiers de $\mathbb K$ et $a$ un élément de 
$\mathcal U_\mathbb K^+$. On considère l'application $m_a:x\mapsto a\cdot x$
qui est un isomorphisme du $\Q$-espace vectoriel $\mathbb K$, et on
note $M_a=(m_{i,j})$ sa matrice dans la base $\mathcal B$. Comme $a$ est
un entier de $\mathbb K$, le $\Z$-module $\mathcal O_\mathbb K$ est stable 
par $m_a$,
donc la matrice $M_a$ est entière, et le déterminant de $m_a$ est
par définition égal à la norme de $a$ qui vaut 1. La matrice $M_a$
appartient donc à $\SL_k(\Z)$. Par ailleurs, l'égalité $m_a(b_i)=a\cdot b_i$
peut s'écrire $\sum_jm_{j,i}b_j=a\cdot b_i$ pour tout $i$, c'est-à-dire
que
\begin{equation}
\trans{M_a}\left(\begin{array}{c}b_1\\\vdots\\b_k\end{array}\right)
=a\cdot\left(\begin{array}{c}b_1\\\vdots\\b_k\end{array}\right).
\end{equation}
Si on pose $v=(b_1,\ldots,b_k)$, c'est un vecteur difficilement approchable
car ses composantes forment une base d'un corps de nombres, et c'est
un vecteur propre de la matrice $\trans{M_a}$, pour tout $a\in
\mathcal U_\mathbb K^+$. L'application $\rho:a\mapsto\trans{M_a}$,
qui est naturellement un morphisme de groupe, est donc à valeur dans
$G_v$, et elle est injective car le fait que $M_a=I$ implique que $a=1$.
On peut se convaincre que c'est une bijection en remarquant que 
le morphisme injectif $\varphi:G_v\to\R_+^*$ qui apparaît dans la 
démonstration 
du théorème~\ref{ap:th2} est en fait à valeur dans $\mathcal U_\mathbb K^+$
et constitue une application réciproque de $\rho$.

Considérons maintenant un vecteur $v=(v_1,\ldots,v_k)\in\R^k$ dont les 
composantes sont linéairement indépendantes sur $\Q$. L'espace vectoriel 
$E_v$ n'est pas nécessairement un corps, mais il contient le $\Q$-espace 
vectoriel engendré par l'image du morphisme $\varphi:G_v\to\R_+^*$,
qui est un corps et qu'on notera $\mathbb K$. Si $M=(m_{i,j})\in G_v$ et
$Mv=\lambda v$, alors $E_v$ est stable par multiplication par $\lambda$:
en effet, on a alors $\lambda v_i=\sum_jm_{i,j}v_j\in E_v$, pour tout $i$.
Cette multiplication définit une action libre de $\varphi(G_v)$ sur $E_v$,
qui s'étend par $\Q$-linéarité en une opération multiplicative de 
$\mathbb K$ sur $E_v$, c'est-à-dire que 
$E_v$ est un espace vectoriel sur $\mathbb K$. En notant $d$ le degré de
$\mathbb K$, on se donne une base $(e_1,\ldots,e_{k/d})$ de
$E_v$ sur $\mathbb K$ telle que $e_1=1$. Si $\mathcal B=(b_1,\ldots,
b_d)$ est une base de $\mathbb K$, les produits $(b_i\cdot e_j)_{i,j}$ forment
une base $\mathcal B'$ de $E_v$ sur $\Q$.

L'équation $\lambda v_i=\sum_jm_{i,j}v_j$ montre que $\trans{M}$ est
la matrice de l'endomorphisme $m_\lambda:x\mapsto\lambda\cdot x$ dans la base 
$(v_i)$ de $E_v$; on va chercher quelle forme prend sa matrice dans la base
$\mathcal B'$. Remarquons d'abord que $\mathbb K$ est stable par
$m_\lambda$, ainsi que tous les sous-espaces de la forme $\mathbb K\cdot
e_j$, et que les matrices des restrictions $m_{\lambda|\mathbb K\cdot
e_j}$ dans les bases $(b_i\cdot e_j)_i$ sont toutes les mêmes. Il
suffit donc de déterminer la matrice de $m_\lambda$ restreint à 
$\mathbb K$. Or on l'a déjà calculé, c'est la transposée de la 
matrice $A=\rho(\lambda)\in G_{\mathbb K,\mathcal B}$.
La transposée de la matrice de $m_\lambda$ dans la base $\mathcal B'$
est donc
$$A'=\left(\begin{array}{ccc}A&0&0\\0&\ddots&0\\0&0&A\end{array}
\right).$$
Si on note $P$ la transposée de la matrice de passage de $(v_i)$
à $\mathcal B'$, on a donc $M=P^{-1}A'P$. Comme la matrice $P$ ne dépend pas
de $M$, on a montré que $G_v$ est contenu dans le groupe $P^{-1}G'P$.

Réciproquement, un élément $M=(m_{i,j})\in P^{-1}G'P$ est la transposée de 
la matrice dans la base $(v_i)$ d'un endomorphisme $m_\lambda$
avec $\lambda\in\mathcal U_\mathbb K^+$, et donc vérifie
$\lambda v_i=\sum_jm_{i,j}v_j$ pour tout $i$. Le vecteur $v$ est
donc vecteur propre de $M$, et si de plus $M\in\SL_k(\Z)$, alors
$M$ appartient à $G_v$. Finalement, on a bien $G_v=P^{-1}G'P\cap\SL_k(\Z)$.
\end{demo}

\noindent Pierre \textsc{Jammes}\\
Université d'Avignon\\
laboratoire de mathématiques\\
33 rue Louis Pasteur\\
F-84000 Avignon\\
\texttt{Pierre.Jammes@ens-lyon.org}
\end{document}